\documentclass[12pt]{article}

\usepackage[margin=1in]{geometry}  
\usepackage{graphicx}              
\usepackage{amsmath}               
\usepackage{amsfonts}              
\usepackage{amsthm}                
\usepackage{amssymb}

\newcommand{\Lu}{\mathfrak{u}}
\newcommand{\Lso}{\mathfrak{so}}

\newcommand{\Lc}{\mathfrak{c}}
\newcommand{\Lgl}{\mathfrak{gl}}
\newcommand{\Lsl}{\mathfrak{sl}}
\newcommand{\Lg}{\mathfrak{g}}
\newcommand{\Lh}{\mathfrak{h}}
\begin{document}

\nocite{*}

\title{Endomorphism algebras of Kuga-Satake varieties.}

\author{Evgeny Mayanskiy}

\maketitle

\begin{abstract}
We compute endomorphism algebras of Kuga-Satake varieties associated to $K3$ surfaces.
\end{abstract}



\section{Preliminary remarks.} 

Let $V$ be a $\mathbb Q$-lattice of transcendental cycles on a $K3$ surface $X$, $\phi\colon V{\otimes}_{\mathbb Q} V\rightarrow \mathbb Q$ the polarization of the weight $2$ Hodge structure on $V$, $E=End_{Hdg}(V)$, $\Phi\colon V{\otimes}_{E} V\rightarrow E$ the hermitian or bilinear form constructed in \cite{Zarhin}, $\phi=tr \circ \Phi$.\\

Let $C^{}(V)$ be the Clifford algebra of the quadratic space $(V,\phi)$ over $\mathbb Q$, $C^{+}(V)$ the even Clifford algebra and $KS(X)$ the Kuga-Satake variety of $X$. Here we define $KS(X)$ from the weight $2$ Hodge structure on the lattice of transcendental cycles $V$ rather than on the whole lattice of primitive cycles $H^2(X,\mathbb Q)_{prim}$. In particular, the Kuga-Satake variety defined here is isogenous to a power of the Kuga-Satake variety defined using the whole lattice of primitive cycles (see \cite{KugaSatake}, \cite{Morrison}, \S 4).\\ 

We want to compute the endomorphism algebra $End(KS(X))_{\mathbb Q}=End_{Hdg}(C^{+}(V))$.\\ 

Let $Z(\Phi)$ be the $\mathbb Q$-algebraic group $Res_{E/\mathbb Q}(SO(V,\Phi))$, if $E$ is a totally real field, or $Res_{E_0/\mathbb Q}(U(V,\Phi))$, if $E=E_0(\theta)$ is a CM-field (with the totally real subfield $E_0$). Recall, that according to \cite{Zarhin}, $Z(\Phi)$ is the Hodge group of the Hodge structure on $V$.\\

Let $CSpin(\phi)\colon = \{ g\in C^{+}(V)^{*}  \; \mid \; gVg^{-1}\subset V \}$. Consider the vector representation $\rho \colon CSpin(\phi)\rightarrow GL(V)$, $g\mapsto (v\mapsto gvg^{-1})$ and the spin representation $\sigma\colon CSpin(\phi)\rightarrow GL(C^{+}(V))$, $g\mapsto (x\mapsto gx)$. Let $ZSpin(\Phi)\colon =\{ g\in CSpin(\phi)  \;\mid\; {\rho}(g)\in Z(\Phi)  \}=$ ${\rho}^{-1}(Z(\Phi))\subset CSpin(\phi)$. Note that ${\rho}(ZSpin(\Phi))=Z(\Phi)$.\\

{\bf Lemma 1.} {\it The Mumford-Tate group of the weight $1$ Hodge structure on $C^{+}(V)$ is the preimage with respect to $\rho$ of the Mumford-Tate group of the weight $2$ Hodge structure on $V$.}\\

{\it Proof:} The same as Proposition 6.3 in \cite{vanGeemen}. If $h_X\colon S^{1}\rightarrow GL(V)$ and $h_{KS(X)}\colon S^{1}\rightarrow GL(C^{+}(V))$ denote the corresponding Hodge structures, then $h_X={\rho}\circ {\sigma}^{-1}\circ h_{KS(X)}$ (as shown in \cite{vanGeemen}). {\it QED}\\

{\bf Corollary.} {\it $End(KS(X))_{\mathbb Q}\cong End_{ZSpin(\Phi)}(C^{+}(V))$, where $ZSpin(\Phi)$ acts on $C^{+}(V)$ via the spin representation ${\sigma}{\mid}_{ZSpin(\Phi)}$.}\\

So, if $C^{+}(V)=\bigoplus_{j} T_{j}^{\oplus m_j}$ is the decomposition of ${\sigma}{\mid}_{ZSpin(\Phi)}$ into a direct sum of irreducible (mutually non-isomorphic) representations $T_j$, then $End(KS(X))_{\mathbb Q}\cong \prod_{j} Mat_{m_j\times m_j}(D_j)$ as $\mathbb Q$-algebras, where $D_j=End_{CSpin(\Phi)}(T_j)$.\\

Let us assume that $m=dim_EV\geq 3$, if $E=E_0$ is totally real, and $m=dim_EV \geq 2$, if $E=E_0(\theta)$ is a CM-field. In the totally real case condition $m\geq 3$ is automatically satisfied for any $K3$ surface $X$ (see \cite{Mayanskiy} and \cite{vanGeemen2}). In what follows we will often denote the field of rational numbers $\mathbb Q$ by $k$ and $E_0$ by $L$. Our approach is not invariant in the sense that we choose a basis in $V$ which diagonalizes $\Phi$ right from the start (see Section 2).\\

Consider the epimorphism $\pi\colon CSpin(\phi)\rightarrow SO(\phi)$ of algebraic groups over $\mathbb Q$ (induced by the vector representation $\rho$ above) with fiber $ker(\pi)={\mathbb G}_m\subset CSpin(\phi)$ and its restriction ${\pi}_0\colon Spin(\phi)\rightarrow SO(\phi)$ to the subgroup $Spin(\phi)\subset CSpin(\phi)$. Then ${\pi}_0$ is a double etale covering \cite{BConrad}.\\

The argument above shows that the Hodge group $Hdg$ of the Kuga-Satake structure on $C^{+}(V)$ satisfies inclusions:
$$
Hdg\subset ({{\pi}_0}^{-1}(Z(\Phi)))^{0} \cdot {\mathbb G}_m \; \; \mbox{and} \; \; \; ({{\pi}_0}^{-1}(Z(\Phi)))^{0} \subset Hdg
$$
(hereafter for an algebraic group $G$ we let $G^{0}$ denote the connected component of the identity and $Lie(G)$ the Lie algebra of $G$).\\

Hence the $\mathbb Q$-algebra
$$
End_{Hdg}(C^{+}(V))=End_{({{\pi}_0}^{-1}(Z(\Phi)))^{0}}(C^{+}(V))=End_{Lie({{\pi}_0}^{-1}(Z(\Phi)))}(C^{+}(V))=End_{Lie(Z(\Phi))}(C^{+}(V)).
$$

Let $\Lg=Lie(Z(\Phi))$. Then $\Lg=Res_{E/k}(\Lso(\Phi))$, if $E$ is totally real, or $\Lg=Res_{E_0/k}(\Lu(\Phi))$, if $E=E_0(\theta)$ is a CM-field (${\theta}^2\in E_0$), where $k=\mathbb Q$.\\

Hence what we are looking for is the algebra of intertwining operators $End_{\Lg}(C^{+}(V))$ of the $\mathbb Q$-linear representation of the Lie algebra $\Lg$ over $\mathbb Q$ induced by the spin representation of $\Lso(\phi)$ in $C^{+}(V)$ via the inclusion of Lie algebras $\Lg\subset \Lso(\phi)$ corresponding to the inclusion of the $\mathbb Q$-algebraic groups $Z(\Phi)\subset SO(\phi)$ above.\\

The problem of computing endomorphism algebras of Kuga-Satake varieties was addressed earlier by Bert van Geemen in papers \cite{vanGeemen1} and \cite{vanGeemen2}. In particular, in \cite{vanGeemen1} he considered the case of the CM-field, which is quadratic over $\mathbb Q$ and in \cite{vanGeemen2} he considered the case of the totally real field, computed the endomorphism algebra in several special cases and made some general remarks. A different computation of the endomorphism algebra of the Kuga-Satake variety in the totally real case was done by  Ulrich Schlickewei \cite{Schlickewei}.\\

Our solution uses the same ideas as (some of the ideas) in papers \cite{vanGeemen1} and \cite{vanGeemen2}. We compute the decomposition of the restriction to $\Lg$ of the spin representation of $\Lso(\phi)$ into irreducible subrepresentations over a splitting field of $\Lg$, and then apply Galois descent.\\

Our main result is Theorem 1 in Section 4 complemented by the computation of primary representations (which are the multiples of irreducible representations $T_j$ above) and division algebras (which are the endomorphism algebras of $T_j$) in subsequent sections. In this text a 'primary representation' means a multiple of an irreducible representation. Some general observations regarding representations over arbitrary fields are collected in the Section 2. In Section 3 we introduce Galois-invariant Cartan subalgebras. In Section 4 we compute decompositions of representations over a splitting field. In Section 5 we construct primary representations over $\mathbb Q$ whose irreducible components appear in Theorem 1. In Section 6 we compute the division algebras which are the endomorphism algebras of those irreducible components. Section 7 is devoted to examples.\\

\section{Some remarks on Galois theory of representations.} 

Let $F/k=\mathbb Q$ be a finite Galois extension, $\Lg=\Lc\oplus {\Lg}'$ be a reductive Lie algebra over $k$, $\Lc\subset \Lg$ be its center and ${\Lg}'\subset \Lg$ be its derived subalgebra. Let $S=Gal(F/k)$ and $\Lh\subset \Lg\otimes_k F$ be a Galois-invariant (i.e. such that $g(\Lh)=\Lh$ for any $g\in S$) splitting Cartan subalgebra. Let $B$ be a basis of the root system $R$ of $(\Lg\otimes_k F,\Lh)$. In what follows we assume that all the representations of $\Lg$ we are dealing with are finite-dimensional and can be integrated to representations of a reductive algebraic group with Lie algebra $\Lg$ (in order to guarantee their complete reducibility).\\

Let $\rho\colon \Lg\rightarrow End_k(W)$ be a representation of $\Lg$ over $k$ and $W\otimes_k F=\oplus_{\alpha} V_{\alpha}$ its decomposition into irreducible subrepresentations over $F$. Let ${\rho}_{\alpha}={\rho}{\mid}_{V_{\alpha}}$ be an irreducible representation of $\Lg\otimes_k F$ with primitive element $v_{\alpha}\in W\otimes_k F$ with highest weight ${\omega}_{\alpha}\in Hom_F(\Lh,F)$ (with resprect to $B$). Then for any $g\in S$, ${\rho}_{\alpha}^{g}\colon =\rho {\mid}_{g(V_{\alpha})}$ is an irreducible representation of $\Lg\otimes_k F$ with primitive element $g(v_{\alpha})\in W\otimes_k F$ with highest weight $g\circ {{\omega}_{\alpha}}\circ g^{-1}\in Hom_F(\Lh,F)$ with respect to the basis $g\circ B\circ g^{-1}$ of $R$. Since the Weyl group ${\mathcal W}_{R}$ of $R$ acts simply transitively on the set of bases of $R$, for any $g\in S$ there exists unique $w(g)\in {\mathcal W}_R$ such that $g\circ B\circ g^{-1}=w(g)(B)$. Hence ${\rho}_{\alpha}^{g}$ is an irreducible representation of $\Lg\otimes_k F$ with primitive element $g(v_{\alpha})\in W\otimes_k F$ with highest weight ${\omega}_{\alpha}^{g}\colon = w(g)^{-1}(g\circ {\omega}_{\alpha}\circ g^{-1})\in Hom_F(\Lh,F)$ (with respect to $B$).\\

{\bf Lemma 3.} {\it Suppose that ${\rho}_1\colon \Lg\rightarrow End_k(W_1)$ and ${\rho}_2\colon \Lg\rightarrow End_k(W_2)$ are two irreducible representations of $\Lg$ over $k$, $V_{\alpha}\subset W_1\otimes_k F$ and $V_{\beta}\subset W_2\otimes_k F$ are two irreducible subrepresentations of $\Lg\otimes_k F$ over $F$. Then $W_1\cong W_2$ as $\Lg$-modules over $k$, if and only if there exist ${\sigma}, {\tau}\in S$ such that $({\rho}_1 {\mid}_{V_{\alpha}})^{\sigma}\cong ({\rho}_2 {\mid}_{V_{\beta}})^{\tau}$ as $\Lg\otimes_k F$-modules over $F$.}\\

{\it Proof:} Schur's lemma. {\it QED}\\

{\bf Corollary.} {\it If ${\rho}_{\alpha}\colon \Lg\otimes_k F \rightarrow End_F(V_{\alpha})$ is an irreducible representation of $\Lg\otimes_k F$ over $F$, then there exists at most one irreducible representation $\rho\colon \Lg\rightarrow End_k(W)$ of $\Lg$ over $k$, such that ${\rho}_{\alpha}$ is a subrepresentation of $\rho \otimes_k F$.}\\

Using the notation of the remark preceeding Lemma 3, let $W=\oplus_{\gamma} W_{\gamma}$ be a decomposition of $\rho$ into irreducible subrepresentations over $k$. Then for any $\gamma$ such that $V_{\alpha}\subset W_{\gamma}\otimes_k F$, by Galois descent we have: 
$$
\bigoplus_{{\gamma}'\colon W_{{\gamma}'}\cong W_{{\gamma}}\; \mbox{as}\; \Lg\mbox{-modules} } W_{{\gamma}'}=\left( \bigoplus_{{\alpha}' \colon {\rho}_{\alpha}^{\tau}\cong {\rho}_{{\alpha}'}^{\sigma}\; \mbox{for some} \; \tau, \sigma \in S}  V_{{\alpha}'} \right)^{S}.
$$

Hence $dim_k \left( \bigoplus_{{\gamma}'\colon W_{{\gamma}'}\cong W_{{\gamma}} \;\mbox{as}\; \Lg\mbox{-modules} } W_{{\gamma}'} \right)=$ $dim_F \left( \bigoplus_{{\alpha}'\colon {\rho}_{\alpha}^{\tau}\cong {\rho}_{{\alpha}'}^{\sigma}\; \mbox{for some}\; \tau, \sigma \in S } V_{{\alpha}'} \right)=$ $\sum_{{\alpha}'\colon \exists \tau, \sigma \in S \colon {\omega}_{\alpha}^{\tau}={\omega}_{{\alpha}'}^{\sigma}} dim_F (V_{{\alpha}'})=$ $m_{\alpha}\cdot dim_k(W_{\gamma})$, where $m_{\alpha}$ is the multiplicity of $W_{\gamma}$ in the decomposition above.\\

So, if $W_{{\gamma}_1}, ..., W_{{\gamma}_p}$ are pairwise nonisomorphic (as $\Lg$-modules) irreducible $\Lg$-submodules of $W$ over $k$ (with the corresponding $\Lg\otimes_k F$-submodules $V_{{\alpha}_i}\subset W\otimes_k F$) appearing in the decomposition above, then $W=\oplus_{i} W_{{\gamma}_i}^{\oplus m_{{\alpha}_i}}$ and
$$
End_{\Lg}(W)\cong \prod_i Mat_{m_{{\alpha}_i}\times m_{{\alpha}_i}}(D_i) \; \; \mbox{as}\; k-\mbox{algebras},
$$
where $D_i=End_{\Lg}(W_{{\gamma}_i})$, $W_{{\gamma}_i}$ is the unique irreducible $\Lg$-module over $k$ such that $W_{{\gamma}_i}\otimes_k F$ contains $V_{{\alpha}_i}$ as a $\Lg\otimes_k F$-submodule over $F$ and $m_{{\alpha}_i}= \left(  \sum_{{\alpha}'\colon \exists \sigma \in S \colon {\omega}_{{\alpha}'}={\omega}_{{\alpha}_i}^{\sigma}} dim_F(V_{{\alpha}'}) \right)   / dim_k(W_{{\gamma}_i})$. We can also write:
$$
m_{{\alpha}_i}=\frac{dim_F(V_{{\alpha}_i}) \cdot \sum_{\sigma \in S} mult({\omega}_{{\alpha}_i}^{\sigma}) }{n_{{\omega}_{{\alpha}_i}} \cdot dim_k(W_{{\gamma}_i})},
$$
where $milt(\omega)$ is the multiplicity of the irreducible representation of $\Lg\otimes_k \mathbb C$ with highest weight $\omega$ (relative to the chosen $\Lh$ and $B$) in $W\otimes_k \mathbb C$ and $n_{\omega}$ is the stabilizer of $\omega$ under the action of the Galois group $S=Gal(F/k)$ on weights. Note that $\{ {\omega}_{{\alpha}_i} \}$ is a set of representatives of the orbits of the action of $S$ on the set of highest weights of irreducible representations of $\Lg\otimes_k \mathbb C$ appearing as irreducible components of $W\otimes_k \mathbb C$.\\

This reduces the study of $End_{\Lg}(W)$ to the study of the (uniquely determined) ($k=\mathbb Q$)-forms of irreducible $\Lg\otimes_k \mathbb C$-submodules of $W\otimes_k \mathbb C$ (i.e. $D_i=End_{\Lg}(W_{{\gamma}_i})$ and $dim_k(W_{{\gamma}_i})$) and the description of the Galois action (of the finite group $Gal(F/k)$) on the weights of $\Lg\otimes_k \mathbb C$ over $\mathbb C$.\\

\section{Description of the Galois action, Cartan subalgebras and bases of the root systems.} 

According to Section 2, we need to specify a splitting field $F$ of $\Lg$ (which should be a Galois extension of $k$), a Galois-invariant splitting Cartan subalgebra $\Lh\subset \Lg \otimes_{k} F$ (i.e. $\Lh$ should be $Gal(F/k)$-stable) and a basis $B$ of the root system $R$ of the split reductive Lie algebra $(\Lg \otimes_{k} F, \Lh)$.\\

Let us assume that $\Phi=d_1\cdot X_1^2+...+d_m\cdot X_m^2$ (if $E=E_0=L$ is totally real) or $\Phi=d_1\cdot X_1\bar{X_1}+...+d_m\cdot X_m\bar{X_m}$ (if $E=E_0(\theta)$, ${\theta}^2\in E_0=L$ is a CM-field), where $d_i\in L$ for any $i$. In other words, we reduce the Hermitian (or quadratic) form $\Phi$ to a diagonal form, i.e. choose an orthogonal (with respect to $\Phi$) basis of $V$ such that $X_i$ are the corresponding coordinates.\\

Let $k=\mathbb Q$ and $F / k$ be a finite Galois extension such that $F$ contains $L$, $\sqrt{d_i}$ for any $i$, $\sqrt{-1}$ and $\theta$ (if $E=E_0(\theta)$ is a CM-field, ${\theta}^2\in E_0$).\\

Let $r=[L\colon k]$ and ${\sigma}_1,..., {\sigma}_r \colon L\hookrightarrow F$ be the list of all field embeddings of $L$ into $F$.\\

\subsection{Case of the totally real field.}

Let us consider first the case $\Lg=Res_{L/k}(\Lso(\Phi))\subset \Lso(\phi)$ (i.e. $E=E_0$ is totally real). We will denote by $E_{i,j}$ a matrix with all entries equal to $0$ except for the entry $(i,j)$ which is equal to $1$.\\

Let ${\Lh}_0=Span_L(A_1,...,A_l)$, where $l=[\frac{m}{2}]$ and $A_i=d_{m-i+1}\cdot E_{m-i+1,i}-d_i\cdot E_{i,m-i+1}$, $1\leq i\leq l$. Let ${\Lh}_i={\Lh}_0\otimes_{L,{\sigma}_i} F\subset \Lso(\Phi)\otimes_{L,{\sigma}_i} F$ and $\Lh={\Lh}_1\times ...\times {\Lh}_r\subset \oplus_{i=1}^{r}(\Lso(\Phi)\otimes_{L,{\sigma}_i} F)\cong Res_{L/k}(\Lso(\Phi))\otimes_{k} F=g \otimes_{k} F$. Then $\Lh\subset \Lg \otimes_{k} F$ is a splitting Cartan subalgebra.\\

Note that over $F$ we have $\Phi=d_1\cdot X_1^2+...+d_m\cdot X_m^2=\sum_{i=1}^{l}Y_i\cdot Y_{-i}+{\epsilon} {Y_0}^2$, where $\epsilon =0$, if $m$ is even, $\epsilon =1$, if $m$ is odd, $Y_i=\sqrt{d_i}\cdot X_i+\sqrt{-d_{m-i+1}}\cdot X_{m-i+1}$, $Y_{-i}=\sqrt{d_i}\cdot X_i-\sqrt{-d_{m-i+1}}\cdot X_{m-i+1}$ and $Y_0=\sqrt{d_{l+1}}\cdot X_{l+1}$.\\

This implies that for any $i, j$ we have $A_j\otimes_{L,{\sigma}_i} 1={\Gamma}_j\cdot H_j$, where ${\Gamma}_j=-\sqrt{{\sigma}_i(d_j)}\cdot \sqrt{-{\sigma}_i(d_{m-j+1})}\in F$ ($1\leq j\leq l$) (in future we will be writing $d_j$ instead of ${\sigma}_i(d_j)$) and $H_j=E_{j,j}-E_{-j,-j}$ (using notation form \cite{Bourbaki}, \S 13). Hence for any $i$ subalgebra ${\Lh}_i\subset \Lso(\Phi)\otimes_{L,{\sigma}_i} F$ is the same splitting Cartan subalgebra as in \cite{Bourbaki}, \S 13. By construction $\Lh\subset \Lg\otimes_{k} F$ is Galois-invariant.\\

Let $R_0$ be the root system of type $B_l$, if $m=2l+1$ (respectively, of type $D_l$, if $m=2l$) from \cite{Bourbaki}, \S 13, i.e. $R_0=\{\pm {\epsilon}_p, \pm {\epsilon}_p \pm {\epsilon}_q\}$ (respectively, $R_0=\{\pm {\epsilon}_p \pm {\epsilon}_q\}$) with basis $B_0=\{ {\epsilon}_1-{\epsilon}_2, {\epsilon}_2-{\epsilon}_3,..., {\epsilon}_{l-1}-{\epsilon}_l, {\epsilon}_l \}$ (respectively, $B_0=\{ {\epsilon}_1-{\epsilon}_2, {\epsilon}_2-{\epsilon}_3,..., {\epsilon}_{l-1}-{\epsilon}_l, {\epsilon}_{l-1}+{\epsilon}_l \}$) (using notation from \cite{Bourbaki}, \S 13).\\

Then for any $i$ the root system of $(so(\Phi)\otimes_{L,{\sigma}_i} F, h_0\otimes_{L,{\sigma}_i} F)$ is $R_i=\{\pm {\epsilon}_p\otimes_{L,{\sigma}_i} {\Gamma}_p, \pm {\epsilon}_p\otimes_{L,{\sigma}_i} {\Gamma}_p \pm {\epsilon}_q\otimes_{L,{\sigma}_i} {\Gamma}_q\}$ with basis
$$
B_i=\{ {\epsilon}_1\otimes_{L,{\sigma}_i} {\Gamma}_1-{\epsilon}_2\otimes_{L,{\sigma}_i} {\Gamma}_2, {\epsilon}_2\otimes_{L,{\sigma}_i} {\Gamma}_2-{\epsilon}_3\otimes_{L,{\sigma}_i} {\Gamma}_3,..., {\epsilon}_{l-1}\otimes_{L,{\sigma}_i} {\Gamma}_{l-1}-{\epsilon}_l\otimes_{L,{\sigma}_i} {\Gamma}_l, {\epsilon}_l\otimes_{L,{\sigma}_i} {\Gamma}_l \}
$$
(respectively, $R_i=\{\pm {\epsilon}_p\otimes_{L,{\sigma}_i} {\Gamma}_p \pm {\epsilon}_q\otimes_{L,{\sigma}_i} {\Gamma}_q\}$ with basis
\begin{multline*}
B_i=\{ {\epsilon}_1\otimes_{L,{\sigma}_i} {\Gamma}_1-{\epsilon}_2\otimes_{L,{\sigma}_i} {\Gamma}_2, {\epsilon}_2\otimes_{L,{\sigma}_i} {\Gamma}_2-{\epsilon}_3\otimes_{L,{\sigma}_i} {\Gamma}_3,..., {\epsilon}_{l-1}\otimes_{L,{\sigma}_i} {\Gamma}_{l-1}-{\epsilon}_l\otimes_{L,{\sigma}_i} {\Gamma}_l,\\
{\epsilon}_{l-1}\otimes_{L,{\sigma}_i} {\Gamma}_{l-1}+ {\epsilon}_l\otimes_{L,{\sigma}_i} {\Gamma}_l \}).
\end{multline*}
Then $R=R_1\sqcup ... \sqcup R_r$ is the root system of $(\Lg\otimes_k F, \Lh)$ and as a basis we can take $B=B_1\sqcup ... \sqcup B_r\subset R$.\\

The action of the Galois group $S=Gal(F/k)$ on weights reduces to its action by permutation on factors of $R_1\times ... \times R_r$ (or on the left cosets $Gal(F/k)/Gal(F/L)$) and to switching signes in front of various ${\Gamma}_p$.\\

Note the isomorphism of root systems $R\cong R_0\sqcup ... \sqcup R_0$ ($r$ factors) under which basis $B$ is identified with $B_0\sqcup ... \sqcup B_0$ ($r$ factors).\\

Let $w_p\in {\mathcal W}_{R_0}$ (where ${\mathcal W}_{R}$ denotes the Weyl group of a root system $R$) be the element of the Weyl group such that $w_p(B_0)={\sigma}_p(B_0)$, where ${\sigma}_p$ is a linear transformation of the $\mathbb Q$-vector space generated by the roots of $R_0$ which switches the sign in front of ${\epsilon}_p$ and does not change other ${\epsilon}_q$'s. Then in the notation of Section 2 for any $g\in S$, ${\omega}_{\alpha}^{g}=(\prod_{p\in P_1(g)} {w_p}^{-1})\sqcup ... \sqcup (\prod_{p\in P_r(g)} {w_p}^{-1}) (g\circ {\omega}_{\alpha} \circ g^{-1})\in Hom_F(\Lh,F)$, where $P_i(g)=\{ p \; \mid \; g^{-1}({\epsilon}_p \otimes_{L,{\sigma}_i} {\Gamma}_p)=-{\epsilon}_p \otimes_{L,g^{-1}\circ {\sigma}_i} {\Gamma}_p \}$.\\

\subsection{Case of the CM-field.}

Now let us consider the case $\Lg=Res_{L/k}(\Lu(\Phi))\subset \Lso(\phi)$ (i.e. $E=E_0(\theta)$ is a CM-field, ${\theta}^2\in E_0=L$).\\

Let ${\Lh}_0=Span_L(A_1,...,A_m)$, where $A_i=\theta\cdot E_{i,i}$, ${\Lh}_i={\Lh}_0\otimes_{L,{\sigma}_i} F \subset \Lu(\Phi)\otimes_{L,{\sigma}_i} F\cong \Lgl(m, F)$ and $\Lh={\Lh}_1\times ... \times {\Lh}_r\subset \oplus_{i=1}^{r} (\Lu(\Phi)\otimes_{L,{\sigma}_i} F)\cong \Lgl(m,F)^{\otimes r}\cong Res_{L/k}(\Lu(\Phi))\otimes_k F=\Lg \otimes_k F$. Then $\Lh\subset \Lg \otimes_k F$ is a splitting Cartan subalgebra.\\

Note that over $F$ we have $\Phi=d_1\cdot X_1\bar{X_1}+...+d_m\cdot X_m\bar{X_m}=Y_1\bar{Y_1}+...+Y_m\bar{Y_m}$, where $Y_i=\sqrt{d_i}\cdot X_i$. Hence for any $i, j$ we have $A_j\otimes_{L,{\sigma}_i} 1=\theta \cdot E_{j,j}$ (more precisely we have to write $\sqrt{{{\sigma}_i}({\theta}^2)}$ instead of $\theta$ here) in $\Lu(\Phi)\otimes_{L,{\sigma}_i} F\cong \Lgl(m,F)$ and so for any $i$ subalgebra ${\Lh}_i\subset \Lu(\Phi)\otimes_{L,{\sigma}_i} F\cong \Lgl(m,F)$ is the same splitting Cartan subalgebra as in \cite{Bourbaki}, \S 13. By construction $\Lh\subset \Lg\otimes_k F$ is Galois-invariant.\\

Let $R_0$ be the root system of type $A_{m-1}$ (for the reductive Lie algebra $\Lgl(m)=\Lc\oplus \Lsl(m)$, where $\Lc\subset \Lgl(m)$ is the center), i.e. $R_0=\{ {\epsilon}_p-{\epsilon}_q \}_{p\neq q}$ with basis $B_0=\{ {\epsilon}_1-{\epsilon}_2, ... ,{\epsilon}_{m-1}-{\epsilon}_m  \}$.\\

Then for any $i$ the root system of $(\Lu(\Phi)\otimes_{L,{\sigma}_i} F, {\Lh}_i)\cong (\Lgl(m),$ diagonal matrices $)$ is $R_i=\{ {\epsilon}_p\otimes_{L,{\sigma}_i} \theta-{\epsilon}_q\otimes_{L,{\sigma}_i} \theta  \}$ with basis $B_i=\{ {\epsilon}_1\otimes_{L,{\sigma}_i} \theta-{\epsilon}_2\otimes_{L,{\sigma}_i} \theta, ... ,{\epsilon}_{m-1}\otimes_{L,{\sigma}_i} \theta-{\epsilon}_m\otimes_{L,{\sigma}_i} \theta  \}$. Then $R=R_1\sqcup ... \sqcup R_r$ is the root system of $(\Lg\otimes_k F, \Lh)$ and as a basis we can take $B=B_1\sqcup ... \sqcup B_r \subset R$.\\

The action of the Galois group $S=Gal(F/k)$ on weights reduces to its action by permutation on factors of $R_1\sqcup ... \sqcup R_r$ (or on the left cosets $Gal(F/k)/Gal(F/L)$) and to multiplication of various $\theta$ by $-1$.\\

Note the isomorphism of root systems $R\cong R_0\sqcup ... \sqcup R_0$ ($r$ factors) under which basis $B$ is identified with $B_0\sqcup ... \sqcup B_0$ ($r$ factors).\\

Let $w_0 \in {\mathcal W}_{R_0}$ be such that $w_0(B_0)=-B_0$. Then in the notation of Section 2 for any $g\in S$, ${\omega}_{\alpha}^{g}={(w_0)^{-P_1(g)}}\sqcup ... \sqcup {(w_0)^{-P_r(g)}} (g\circ {\omega}_{\alpha} \circ g^{-1})\in Hom_F(\Lh,F)$, where $P_i(g)=1$, if $g^{-1}({{\epsilon}_p}\otimes_{L,{\sigma}_i} \theta)=- {{\epsilon}_p}\otimes_{L,g^{-1}\circ{\sigma}_i} \theta$ and $P_i(g)=0$ otherwise (the action of $w_0$ is extended to the center of $\Lgl(m,F)$ as multiplication by $-1$).\\

\section{Decomposition of the restriction of the spin representation over a splitting field.} 

In order to apply the general statements of Section 2, we need to decompose the $F$-linear extension of the restriction of the spin representation of $\Lso(\phi)$ in $C^{+}(V)$ to $\Lg\subset \Lso(\phi)$ over $F$. For this we need to describe the embedding of Cartan subalgebras induced by the embedding of Lie algebras $\Lg\otimes_k F\subset \Lso(\phi)\otimes_k F$.\\

{\bf Lemma 2.} {\it If $E$ is totally real, then the Lie algebra homomorphism $\oplus_{i=1}^{r}\Lso(\Phi)\otimes_{L,{\sigma}_i} F \subset \Lso(\oplus_{i=1}^{r}(\Phi \otimes_{L, {\sigma}_i} F))=\Lso(\phi)\otimes_k F$ sends $(M_1,...,M_r)$ to $diag(M_1,...,M_r)$.\\

If $E=E_0(\theta)$ is a CM-field (and ${\theta}^2\in E_0=L$ as usual), then the Lie algebra homomorphism $\oplus_{i=1}^{r} \Lgl(m,F)\cong   \oplus_{i=1}^{r}\Lu(\Phi)\otimes_{L,{\sigma}_i} F \subset \Lso(\oplus_{i=1}^{r}(   (\Phi \otimes_{E, {\sigma}_i} F)\oplus   (\Phi \otimes_{E, \bar{\sigma}_i} F)   ))=\Lso(\phi)\otimes_k F$ (where in the last formula ${\sigma}_i$ and $\bar{\sigma}_i$ denote the two extensions of ${\sigma}_i$ to an embedding of $E/k$ into $F/k$) sends $(M_1,...,M_r)$ to $diag(M_1, -\Phi\cdot {M_1}^T\cdot {\Phi}^{-1},...,M_r, -\Phi\cdot {M_r}^T\cdot {\Phi}^{-1})$.}\\

{\it Proof:} One should notice that $Res_{L/k}$ on vector spaces over $L$ is the forgetful functor to the vector spaces over $k$. Hence on the $Res_{L/k}(\Lso(\Phi))$ (respectively, $Res_{L/k}(\Lu(\Phi))$), which is the Galois-invariant subspace of the source, our homomorphisms have exactly the form needed. Extending scalars to $F$ gives the result. See also Proposition 3.8 in \cite{vanGeemen1} and \cite{vanGeemen2}. {\it QED}\\

\subsection{Case of the totally real field.}

Let $E=L$ be a totally real field.\\

For any $i=1,...,r$, $j=1,...,l$ (where $l=\left[ \frac{m}{2} \right]$) let $\hat{H}_{j}^{i}={{\sigma}_i(d_{m-j+1})}\cdot E_{m-j+1+m(i-1),j+m(i-1)}-{{\sigma}_i(d_j)}\cdot E_{j+m(i-1),m-j+1+m(i-1)}\in \Lso(\phi)\otimes_k F$. Then $\hat{H}_{j}^{i}$ are linearly independent elements of the splitting Cartan subalgebra $\hat{\Lh}\subset \Lso(\phi)\otimes_k F$ described in \cite{Bourbaki}, \S 13. They form a basis of $\hat{\Lh}$, if $m$ is even or $r=1$. If $m$ is odd and $r\geq 2$, then $\hat{H}_{j}^{i}$ together with $\hat{H}_{l+1}^{1},...,\hat{H}_{l+1}^{[\frac{r}{2}]}$ form a basis of $\hat{\Lh}$, if we take $\hat{H}_{l+1}^{i}={{\sigma}_{r-i+1}(d_{l+1})}\cdot E_{(l+1)(r-i+1),(l+1)i}-{{\sigma}_{i}(d_{l+1})}\cdot E_{(l+1)i,(l+1)(r-i+1)}$, $1\leq i\leq [\frac{r}{2}]$.\\

Let us denote by $\{  \hat{\epsilon}_{j}^{i} \}$ the corresponding dual basis of $\hat{\Lh}^{*}=Hom_F(\hat{\Lh},F)$. Its elements differ from the elements of the corresponding basis of the dual Cartan subalgebra considered in \cite{Bourbaki}, \S 13 by scalar factors of the form $-\sqrt{{{\sigma}_i}(d_j)}\cdot \sqrt{-{{\sigma}_i}(d_{m-j+1})}$.\\

Lemma 2 above implies that the restriction of $\hat{\epsilon}_{l+1}^{i}$ to the Cartan subalgebra $\Lh\subset \Lg\otimes_k F$ is zero, while for any $j\leq l$ the restriction of $\hat{\epsilon}_{j}^{i}$ to $\Lh$ is the corresponding element of the dual basis of ${\Lh}^{*}$ of the basis $\{ A_j\otimes_{L,{\sigma}_i} 1 \; \mid \; 1\leq i\leq r, \; 1\leq j\leq l \}$ of $\Lh$.\\

If $m\cdot r=dim_k(V)\geq 5$, then according to \cite{Bourbaki}, \S 13 the weights of the spin representation of $\Lso(\phi)\otimes_k F$ in $C^{+}(V)\otimes_k F$ ($V$ is considered as a vector space over $k$) are $\frac{1}{2}\sum_{i,j} {\hat{\epsilon}}_{j}^{i} - \sum_{(i,j)\in I} {\hat{\epsilon}}_{j}^{i}$, where $I$ runs over the subsets of the set of parameters $i$ and $j$ (i.e. $I\subset \{ (i,j) \; \mid \; 1\leq j\leq l\; \mbox{and}\; 1\leq i\leq r \; \mbox{or (if }\; m\; \mbox{is odd and }\; r\geq 2)\; j=l+1\;\mbox{ and }\; 1\leq i\leq [\frac{r}{2}] \}$) and each weight has multiplicity $\frac{dim_k(C^{+}(V))}{2^{[{mr}/{2}]}}=2^{mr-1-[\frac{mr}{2}]}$.\\

As it was remarked in \cite{vanGeemen2}, Lemma 5.5, this implies (if $m\geq 5$) that the restrictions of these weights to $h\subset \hat{h}$ are exactly the weights of the exterior tensor product of the spin representations of $\Lso(\Phi)\otimes_{L,{\sigma}_i} F$ in $C^{+}(V)\otimes_{L,{\sigma}_i} F$ ($V$ is considered as a vector space over $L$), $1\leq i\leq r$, taken with multiplicity $\frac{2^{mr-1-[mr/2]}}{(2^{m-1-l})^r}=2^{r-1}$, if $m$ is even, or with multiplicity $\frac{2^{mr-1-[mr/2]}}{(2^{m-1-l})^r}\cdot 2^{[\frac{r}{2}]}=2^{r-1}$, if $m$ is odd.\\

{\bf Corollary 1.} {\it If $E=E_0=L$ is totally real, then the restriction of the spin representation $\rho\colon \Lso(\phi)\otimes_k F\rightarrow End_F(C^{+}(V\otimes_k F))$ to $\Lg\otimes_k F=\oplus_{i=1}^{r}(\Lso(\Phi)\otimes_{L,{\sigma}_i} F)\subset \Lso(\phi)\otimes_k F$ is the exterior tensor product $\Gamma\cdot ({\rho}_1 \boxtimes ... \boxtimes {\rho}_r)$ of spin representations ${\rho}_i \colon \Lso(\Phi)\otimes_{L,{\sigma}_i} F \rightarrow End_F(C^{+}(V \otimes_{L,{\sigma}_i} F))$ with multiplicity $\Gamma = 2^{r-1}$.}\\

\subsection{Case of the CM-field.}

Let $E=E_0(\theta), {\theta}^2\in E_0=L$ be a CM-field.\\

For any $i=1,...,r$, $j=1,...,m$ let $\hat{H}_{j}^{i}=E_{j+2m(i-1),j+2m(i-1)}- E_{j+m+2m(i-1),j+m+2m(i-1)}\in \Lso(\phi)\otimes_k F$. Then $\hat{H}_{j}^{i}$ form a basis of the splitting Cartan subalgebra $\hat{\Lh}\subset \Lso(\phi)\otimes_k F$ described in \cite{Bourbaki}, \S 13. Let us denote by $\{  \hat{\epsilon}_{j}^{i} \}$ the corresponding dual basis of $\hat{\Lh}^{*}=Hom_F(\hat{\Lh},F)$. This is the same Cartan subalgebra and the same basis as considered in \cite{Bourbaki}, \S 13.\\

Lemma 2 above implies that the restriction of $\hat{\epsilon}_{j}^{i}$ to the Cartan subalgebra $\Lh\subset \Lg\otimes_k F$ is the element $(0,...,\hat{\epsilon}_j,...,0)$ (with $0$ outside of the $i$-th spot) of the Cartan subalgebra (consisting of diagonal matrices) of $\Lgl(m,F)^{\oplus r}$, where $\hat{\epsilon}_j\cong E_{j,j}\in \Lgl(m,F)$ is the $j$-th element of the dual basis of the Cartan subalgebra of $\Lgl(m,F)$ considered in \cite{Bourbaki}, \S 13.\\

If $m\cdot r=\frac{1}{2}\cdot dim_k(V)\geq 3$, then according to \cite{Bourbaki}, \S 13 the weights of the spin representation of $\Lso(\phi)\otimes_k F$ in $C^{+}(V)\otimes_k F$ ($V$ is considered as a vector space over $k$) are $\frac{1}{2}\sum_{i,j} {\hat{\epsilon}}_{j}^{i} - \sum_{(i,j)\in I} {\hat{\epsilon}}_{j}^{i}$. Here $I$ runs over the subsets of $[1,...,r]\times [1,...,m]$. Each weight has multiplicity $\frac{dim_k(C^{+}(V))}{2^{mr}}=2^{mr-1}$ (\cite{Bourbaki}, \S 13).\\

Suppose $m\geq 2$. Then the restrictions of these weights to $h\subset \hat{h}$ are exactly the weights of the exterior tensor product of the exterior algebra representations of $\Lu(\Phi)\otimes_{L,{\sigma}_i} F\cong \Lgl(m,F)$ in ${\wedge}_{E}^{*}(V)\otimes_{E,{\sigma}_i} F$ ($V$ is considered as a vector space over $E$) twisted by $D^{-1/2}$, $1\leq i\leq r$. Here $D^{c}$, $c\in\mathbb Q$ denotes the representation of $\Lu(\Phi)\otimes_{L,{\sigma}_i} F\cong \Lgl(m,F)=\Lc\oplus \Lsl(m,F)$ in ${\wedge}_{E}^{m}(V)\otimes_{E,{\sigma}_i} F\cong {\wedge}_{F}^{m}(V\otimes_{E,{\sigma}_i} F)$ such that $\Lsl(m,F)$ acts trivially, while $1\in F\cong \Lc$ acts as $c\cdot Id$. In other words, $D^{c}\colon \Lgl(m,F)\rightarrow End_F({\wedge}_{E}^{m}(V)\otimes_{E,{\sigma}_i} F)$, $M \mapsto c\cdot Tr(M)\cdot Id$.\\

Indeed, for any $i$, $\sum_{j}\hat{\epsilon}_{j}^{i}$ restricts to $0$ to the Cartan subalgebra of the semi-simple part $\Lsl(m,F)\subset\Lgl(m,F)\cong \Lu(\Phi)\otimes_{L,{\sigma}_i} F$ and to $m\cdot Id_F$ to the center $F\cong \Lc\subset \Lgl(m,F)\cong \Lu(\Phi)\otimes_{L,{\sigma}_i} F$.\\

The exterior tensor product above has multiplicity $\Gamma = 2^{mr-1}$. Indeed, $dim_F(C^{+}(V)\otimes_k F)=2^{2mr-1}$ and $dim_F({\wedge}_{E}^{*}(V)\otimes_{E,{\sigma}_i} F)=2^m$. Hence the dimention of the exterior tensor product is $(dim_F({\wedge}_{E}^{*}(V)\otimes_{E,{\sigma}_i} F))^r=2^{mr}$ and so the multiplicity is $2^{2mr-1}/2^{mr}=2^{mr-1}$.\\ 

{\bf Corollary 2.} {\it If $E=E_0(\theta), {\theta}^2\in E_0=L$ is a CM-field, then the restriction of the spin representation $\rho\colon \Lso(\phi)\otimes_k F\rightarrow End_F(C^{+}(V\otimes_k F))$ to $\Lg\otimes_k F=\oplus_{i=1}^{r}(\Lu(\Phi)\otimes_{L,{\sigma}_i} F)\cong {\Lgl(m,F)}^{\oplus r}\subset \Lso(\phi)\otimes_k F$ is the exterior tensor product $\Gamma\cdot ({\rho}_1 \boxtimes ... \boxtimes {\rho}_r)$ of exterior algebra representations ${\rho}_i \colon \Lgl(m,F) \rightarrow End_F({\wedge}_{F}^{*}(V\otimes_{E,{\sigma}_i} F)\otimes_F F )$ twisted by one-dimensional representations $D^{-1/2}\colon \Lgl(m,F)\rightarrow End_F(F)\cong F$, $M\mapsto (-\frac{1}{2m})\cdot Tr(M)$ with multiplicity $\Gamma = 2^{mr-1}$.}\\

{\bf Remark.} ${\rho}_i$ is a double-valued 'spin' representation of $GL(m,F)$.\\

From these Corollaries one can deduce the highest weights of irreducible subrepresentations over $F$ of the restriction to $\Lg\otimes_k F \subset \Lso(\phi)\otimes_k F$ of the spin representation $\rho\colon \Lso(\phi) \rightarrow End_k(C^{+}(V))$. Then one can use the description of the Galois action of $S=Gal(F/k)$ on weights of $\Lg\otimes_k F$ given above in order to break down the highest weights into orbits $\{ S\cdot {\omega}_1,...,S\cdot {\omega}_t  \}$. Let us denote the dimension of the irreducible representation of $\Lg\otimes_k F$ with highest weight ${\omega}_i$ by $d_i$. Let $\hat{\rho}_i \colon \Lg\rightarrow End_k(W_i)$ be the (unique) irreducible representation of $\Lg$ over $k$ such that $W_i\otimes_k F$ contains the irreducible representation of $\Lg\otimes_k F$ with highest weight ${\omega}_i$ as a $(\Lg\otimes_k F)$-submodule. Then our analysis in Section 2 implies:\\

{\bf Theorem 1.} {\it
$$
End(KS(X))_{\mathbb Q} \cong End_{\Lg}(W)\cong \prod_i Mat_{m_i\times m_i}(D_i) \; \; \mbox{as}\; {\mathbb Q}-\mbox{algebras},
$$
where $D_i=End_{\Lg}(W_{i})$, $m_i=(d_i/ dim_k(W_i))\cdot \sum_{\omega\in S\cdot {\omega}_i} mult(\omega)$ and $mult(\omega)$ is the multiplicity of the irreducible subrepresentation of the representation of $\Lg\otimes_k F$ on $C^{+}(V\otimes_k F)$ with highest weight $\omega$}.\\

{\bf Remark. } In the analysis above we assumed that $m=dim_EV\geq 5$ (if $E$ is totally real) or $m\geq 2$ (if $E$ is a CM-field and $r=[E:k]/2\geq 2$) or $m\geq 3$ (if $E$ is a CM-field and $r=[E:k]/2=1$). In the case of small $m$ Lie algebras we consider 'degenerate' and requre a separate consideration.\\

\section{$\mathbb Q$-forms of spin representations.} 

Let us describe more explicitely $\mathbb Q$-forms $W_i$ above or at least the corresponding primary representations. We will use corestriction of algebraic structures, as in \cite{vanGeemen2}, \S 6 and (in the case of totally real fields) representation spaces which we are going to construct in the following subsection.\\

\subsection{Galois-invariant sums of ideals of Clifford algebra.}

Let $k=\mathbb Q$, $E=L$ be a totally real number field, $r=[L:k]$. Let $\Phi=d_1\cdot X_1^2+...+d_m\cdot X_m^2$ with respect to basis $\{ e_1,...,e_m   \}$ of $V$, $m=dim_LV$. Let $F/k$ be a finite Galois extension containing $L$, $\sqrt{-1}$ and $\sqrt{d_i}$ for all $i$. Let ${\sigma}_1,...,{\sigma}_r\colon L\hookrightarrow F$ be all the field embeddings over $k$.\\

Let $f_i=\frac{1}{\sqrt{d_i}}\cdot e_i+\frac{1}{\sqrt{-d_{m-i+1}}}\cdot e_{m-i+1}$, $f_{-i}=\frac{1}{\sqrt{d_i}}\cdot e_i-\frac{1}{\sqrt{-d_{m-i+1}}}\cdot e_{m-i+1}$, $1\leq i\leq l=[\frac{m}{2}]$ and $f_0=\frac{1}{\sqrt{d_{l+1}}}\cdot e_{l+1}$. Then $\{ f_i, \; f_{-i}\; \mid \; 1\leq i\leq l \}$ (if $m$ is even) or $\{ f_0,\; f_i, \; f_{-i}\; \mid \; 1\leq i\leq l \}$ (if $m$ is odd) is a basis of $V\otimes_{L,{\sigma}_i} F$, where we denote ${{\sigma}_i}(d_j)$ by $d_j$. With respect to this basis $\Phi=2\sum_{i=1}^{l}Y_i\cdot Y_{-i}+{\epsilon}Y_0^2$, where ${\epsilon}=(1-(-1)^m)/2$.\\

\subsubsection{Even dimension.}

Assume that $m$ is even. Let $f_{{\alpha}_1,...,{\alpha}_l}^i=f_{{\alpha}_1\cdot 1}\cdot...\cdot f_{{\alpha}_l\cdot l}\in C(V\otimes_{L,{\sigma}_i} F)$ for various ${\alpha}_i\in \{  \pm 1 \}$ and $I_{{\alpha}_1,...,{\alpha}_l}^{i}=C(V\otimes_{L,{\sigma}_i} F)\cdot f_{{\alpha}_1,...,{\alpha}_l}^i$, $1\leq i\leq r$. $I_{{\alpha}_1,...,{\alpha}_l}^{i}$ are left ideals of the Clifford algebra $C(V\otimes_{L,{\sigma}_i} F)$ viewed as $F$-vector subspaces.\\

Consider the direct sum of $F$-vector spaces 
$$
\tilde{C}(V\otimes_{L,{\sigma}_i} F)=\tilde{C}(V)\otimes_{L,{\sigma}_i} F=\bigoplus_{{\alpha}_1,...,{\alpha}_l\in \{  \pm 1 \}} I_{{\alpha}_1,...,{\alpha}_l}^{i}.
$$

Note that $g(f_i)\in \{ \pm f_i, \pm f_{-i}  \}$ for any $i$ and $g\in S$. Hence the Galois group $S=Gal(F/k)$ acts on $\tilde{C}(V\otimes_{L,{\sigma}_i} F)$ (by sending an element of the summand $I_{{\alpha}_1,...,{\alpha}_l}$ to its image under the action of $S$ on $C(V\otimes_{L,{\sigma}} F)$ viewed as an element of the summand $I_{{\beta}_1,...,{\beta}_l}$, where $f_{{\beta}_1,...,{\beta}_l}$ is upto a scalar factor the image of $f_{{\alpha}_1,...,{\alpha}_l}$).\\

It follows from the construction that $F$-vector subspaces $\oplus_{i=1}^{r} I_{{\alpha}_1^i,...,{\alpha}_l^i}^{i}\subset \oplus_{i=1}^r \tilde{C}(V)\otimes_{L,{\sigma}_i} F$ for various choices of ${\alpha}_j^i\in \{ \pm 1 \}$ are permuted among themselves under the action of the Galois group $S=Gal(F/k)$.\\

{\bf Remark.} For any ${\alpha}_1,...,{\alpha}_l$ the left ideal $I_{{\alpha}_1,...,{\alpha}_l}^{i}\subset C(V\otimes_{L,{\sigma}_i} F)$ is an $(\Lso(\Phi)\otimes_{L,{\sigma}_i} F )$-subrepresentation of the spin representation, which is either irreducible (if $m$ is odd) or is the sum of two irreducible and non-isomorphic (semi-spin) representations \cite{Chevalley}, \cite{FultonHarris}. In the latter case, let us write $I_{{\alpha}_1,...,{\alpha}_l}^{i}=I_{{\alpha}_1,...,{\alpha}_l}^{i, +}\oplus I_{{\alpha}_1,...,{\alpha}_l}^{i, -}$ for the corresponding (unique) decomposition.\\ 

\subsubsection{Odd dimension.}

Assume that $m$ is odd. Let $f_{{\alpha}_1,...,{\alpha}_l,\gamma}^i=f_{{\alpha}_1\cdot 1}\cdot...\cdot f_{{\alpha}_l\cdot l} \cdot (1+\gamma\cdot f_0) \in C(V\otimes_{L,{\sigma}_i} F)$ for various ${\alpha}_i,\gamma \in \{  \pm 1 \}$ and $I_{{\alpha}_1,...,{\alpha}_l,\gamma}^{i}=C(V\otimes_{L,{\sigma}_i} F)\cdot f_{{\alpha}_1,...,{\alpha}_l,\gamma}^i$, $1\leq i\leq r$. $I_{{\alpha}_1,...,{\alpha}_l,\gamma}^{i}$ are left ideals of the Clifford algebra $C(V\otimes_{L,{\sigma}_i} F)$ viewed as $F$-vector subspaces.\\

Consider the direct sum of $F$-vector spaces 
$$
\tilde{C}(V\otimes_{L,{\sigma}_i} F)=\tilde{C}(V)\otimes_{L,{\sigma}_i} F=\bigoplus_{{\alpha}_1,...,{\alpha}_l,\gamma \in \{  \pm 1 \}} I_{{\alpha}_1,...,{\alpha}_l,\gamma}^{i}.
$$

Note that $g(1+\gamma \cdot f_0)=(1\pm\gamma \cdot f_0)$ for any $g\in S$. Hence the Galois group $S=Gal(F/k)$ acts on $\tilde{C}(V\otimes_{L,{\sigma}_i} F)$ (by sending an element of the summand $I_{{\alpha}_1,...,{\alpha}_l,\gamma}$ to its image under the action of $S$ on $C(V\otimes_{L,{\sigma}} F)$ viewed as an element of the summand $I_{{\beta}_1,...,{\beta}_l,{\gamma}'}$, where $f_{{\beta}_1,...,{\beta}_l,{\gamma}'}$ is upto a scalar factor the image of $f_{{\alpha}_1,...,{\alpha}_l,\gamma}$).\\

It follows from the construction that $F$-vector subspaces $\oplus_{i=1}^{r} I_{{\alpha}_1^i,...,{\alpha}_l^i,{\gamma}^i}^{i}\subset \oplus_{i=1}^r \tilde{C}(V)\otimes_{L,{\sigma}_i} F$ for various choices of ${\alpha}_j^i, {\gamma}^i\in \{ \pm 1 \}$ are permuted among themselves under the action of the Galois group $S=Gal(F/k)$.\\

{\bf Remark.} For any ${\alpha}_1,...,{\alpha}_l,\gamma$ the left ideal $I_{{\alpha}_1,...,{\alpha}_l,\gamma}^{i}\subset C(V\otimes_{L,{\sigma}_i} F)$ is an irreducible $(\Lso(\Phi)\otimes_{L,{\sigma}_i} F )$-subrepresentation of the spin representation (since $m$ is odd by assumption) \cite{Chevalley}, \cite{FultonHarris}.\\ 

We will use $\tilde{C}(V\otimes_{L,{\sigma}_i} F)$ as representation spaces of $(\Lso(\Phi)\otimes_{L,{\sigma}_i} F )$ (the direct sum of its representations on the left ideals of the Clifford algebra) in order to construct primary $\mathbb Q$-forms of spin representations.\\

\subsection{Case of the totally real field and odd dimension.}

Let $E=E_0=L$ be totally real and $m=dim_LV$ odd. Let ${\Sigma}_i\subset C^{+}(V\otimes_{L,{\sigma}_i} F)$, $1\leq i\leq r$ be the irreducible subrepresentation of the spin representation of $\Lso(\Phi) \otimes_{L,{\sigma}_i} F$. Then ${\Sigma}_1 \otimes_F ... \otimes_F {\Sigma}_r$ is an irreducible representation of $\oplus_{i=1}^{r}(\Lso(\Phi)\otimes_{L,{\sigma}_i} F)=\Lg \otimes_k F$.\\

Let $\tilde{C}(V\otimes_{L,{\sigma}_i} F)=\oplus_p S_{p}^{i}$ be a decomposition into irreducible components of the representation of  $\Lso(\Phi)\otimes_{L,{\sigma}_i} F$ considered above. Let ${\Omega}'$ be the finite set of $F$-vector subspaces of $\tilde{C}(V\otimes_{L,{\sigma}_1} F)\otimes_F ... \otimes_F \tilde{C}(V\otimes_{L,{\sigma}_r} F)$ (or of $C(V\otimes_{L,{\sigma}_1} F)\otimes_F ... \otimes_F C(V\otimes_{L,{\sigma}_r} F)$) of the form ${S}_{p_1}^{1} \otimes_F ... \otimes_F {S}_{p_r}^{r}$ for various $p_1,...,p_r$. These subspaces are irreducible subrepresentations of the exterior tensor product of spin representations as a representation of $\oplus_{i=1}^{r} (\Lso(\Phi)\otimes_{L,{\sigma}_i} F)$.\\

Galois group $S=Gal(F/k)$ acts on ${\Omega}'$. Take any element ${S}_{p_1}^{1} \otimes_F ... \otimes_F {S}_{p_r}^{r}$ of ${\Omega}'$. Let $U\subset \tilde{C}(V\otimes_{L,{\sigma}_1} F)\otimes_F ... \otimes_F \tilde{C}(V\otimes_{L,{\sigma}_r} F)$ be the sum of the elements of ${\Omega}'$ (as subspaces of $\tilde{C}(V\otimes_{L,{\sigma}_1} F)\otimes_F ... \otimes_F \tilde{C}(V\otimes_{L,{\sigma}_r} F)$) lying in the $S$-orbit of ${S}_{p_1}^{1} \otimes_F ... \otimes_F {S}_{p_r}^{r}$. Then $U\subset \tilde{C}(V\otimes_{L,{\sigma}_1} F)\otimes_F ... \otimes_F \tilde{C}(V\otimes_{L,{\sigma}_r} F)$ is an $S$-submodule.\\

Since the actions of $\Lg\subset \Lg\otimes_k F$ and $S=Gal(F/k)$ commute, by Galois descent 
$$
\left( U \right)^S\cong \left(  ({\Sigma}_1 \otimes_F ... \otimes_F {\Sigma}_r)^{\oplus n_0} \right)^{S}
$$
is a primary representation of $\Lg$ over $k$ of dimension $n_0\cdot 2^{l\cdot r}$, which contains ${\Sigma}_1 \otimes_F ... \otimes_F {\Sigma}_r$ after extending scalars to $F$.\\

Multiplicity $n_0$ is the length of the $S$-orbit in ${\Omega}'$ of the chosen element ${S}_{p_1}^{1} \otimes_F ... \otimes_F {S}_{p_r}^{r}$ of ${\Omega}'$.\\

{\bf Remark.} We will use notation introduced above. Consider the action of $S=Gal(F/k)$ on $2^{l+1}$ elements (or more precisely on the lines generated by them) $f_{{\beta}_1,...,{\beta}_l,\gamma}$ of $C(V\otimes_L F)$ for various ${{\beta}_1,...,{\beta}_l,\gamma}$ by sign changes in front of $\sqrt{d_i}$'s and $\sqrt{-d_{m-i+1}}$'s in the definition of $f_i$ in terms of $e_j$ (see notation above). Then (if we choose all $S_{p_i}$ to be the same)  
$$
n_0=\frac{\mbox{order of }\; S=Gal(F/k) }{\mbox{order of the stabilizer of}\; f_{1,...,1,1}}.
$$

\subsection{Case of the totally real field and even dimension.}

Let $E=E_0=L$ be a totally real field and $m=dim_LV$ even. Let ${\Sigma}_i^{+}, {\Sigma}_i^{-} \subset C^{+}(V\otimes_{L,{\sigma}_i} F)$, $1\leq i\leq r$ be irreducible (semi-spin) subrepresentations of the spin representation of $\Lso(\Phi) \otimes_{L,{\sigma}_i} F$.\\

Consider the finite set $\Omega$ of $F$-vector spaces of the form ${\Sigma}_1^{{\alpha}_1} \otimes_F ... \otimes_F {\Sigma}_r^{{\alpha}_r}$ for various ${\alpha}_i \in \{ +,-  \}$. They are exactly the irreducible components of the exterior tensor product of spin representations ${\Sigma}_i={\Sigma}_i^{+}\oplus {\Sigma}_i^{-}\subset C^{+}(V\otimes_{L,{\sigma}_i} F)$ of $\oplus_{i=1}^{r} (\Lso(\Phi)\otimes_{L,{\sigma}_i} F)$ (see \cite{Bourbaki}, \S 13, \cite{Chevalley}, \cite{FultonHarris}). They are also the isomorphism classes of simple $\oplus_{i=1}^{r} (\Lso(\Phi)\otimes_{L,{\sigma}_i} F)$-submodules of $C(V\otimes_{L,{\sigma}_1} F)\otimes_F ... \otimes_F C(V\otimes_{L,{\sigma}_r} F)$. Let $\tilde{C}(V\otimes_{L,{\sigma}_i} F)=\oplus_p S_{p}^{i}$ be a decomposition into irreducible components of the representation of  $\Lso(\Phi)\otimes_{L,{\sigma}_i} F$  considered above. Let ${\Omega}'$ be the finite set of $F$-vector subspaces of $C(V\otimes_{L,{\sigma}_1} F)\otimes_F ... \otimes_F C(V\otimes_{L,{\sigma}_r} F)$ (or of $\tilde{C}(V\otimes_{L,{\sigma}_1} F)\otimes_F ... \otimes_F \tilde{C}(V\otimes_{L,{\sigma}_r} F)$) of the form ${S}_{p_1}^{1} \otimes_F ... \otimes_F {S}_{p_r}^{r}$ for various $p_1,...,p_r$. These subspaces are irreducible subrepresentations of the exterior tensor product of spin representations as a representation of $\oplus_{i=1}^{r} (\Lso(\Phi)\otimes_{L,{\sigma}_i} F)$.\\

Galois group $S=Gal(F/k)$ acts naturally on both $\Omega$ and ${\Omega}'$. Let ${\Omega}_1,...,{\Omega}_u$ be the orbits of $S$ on $\Omega$. For any $i$ choose $({\alpha}_1,...,{\alpha}_r)\in {\Omega}_i$ and define $U_i\subset \tilde{C}(V\otimes_{L,{\sigma}_1} F)\otimes_F ... \otimes_F \tilde{C}(V\otimes_{L,{\sigma}_r} F)$ to be the sum of the elements of ${\Omega}'$ (as subspaces of $\tilde{C}(V\otimes_{L,{\sigma}_1} F)\otimes_F ... \otimes_F \tilde{C}(V\otimes_{L,{\sigma}_r} F)$) lying in the $S$-orbit of any ${S}_{p_1}^{1} \otimes_F ... \otimes_F {S}_{p_r}^{r}$, which is isomorphic to ${\Sigma}_1^{{\alpha}_1} \otimes_F ... \otimes_F {\Sigma}_r^{{\alpha}_r}$ as an $\oplus_{i=1}^{r} (\Lso(\Phi)\otimes_{L,{\sigma}_i} F)$-module.\\

Then $U_i\subset \tilde{C}(V\otimes_{L,{\sigma}_1} F)\otimes_F ... \otimes_F \tilde{C}(V\otimes_{L,{\sigma}_r} F)$ is an $S$-submodule and 
$$
\left( U_i \right)^S\cong \left(  \bigoplus_{({\alpha}_1,...,{\alpha}_r)\in {\Omega}_i} ({\Sigma}_1^{{\alpha}_1} \otimes_F ... \otimes_F {\Sigma}_r^{{\alpha}_r})^{\oplus n_{{\alpha}_1,...,{\alpha}_r}} \right)^{S}
$$
is a primary representation of $\Lg$ over $k$ of dimension $\sum_{({\alpha}_1,...,{\alpha}_r)\in {\Omega}_i} n_{{\alpha}_1,...,{\alpha}_r} \cdot 2^{r\cdot (l-1)}$. These representations $\left( U_i \right)^S$, $1\leq i \leq u$ contain all representations of $\Lg\otimes_k F$ of the form ${\Sigma}_1^{{\alpha}_1} \otimes_F ... \otimes_F {\Sigma}_r^{{\alpha}_r}$ after extending scalars to $F$.\\

Multiplicities $n_{{\alpha}_1,...,{\alpha}_r}$ can be computed as follows:
$$
n_{{\alpha}_1,...,{\alpha}_r}=\frac{\mbox{order of the stabilizer of}\; ({\alpha}_1,...,{\alpha}_r)\in\Omega}{\mbox{order of the stabilizer of}\; ({p}_1,...,{p}_r)\in{\Omega}'}.
$$

{\bf Remark.} We will use notation introduced above. Consider the action of $S=Gal(F/k)$ on $2^l$ elements (or more precisely on the lines generated by them) $f_{{\beta}_1,...,{\beta}_l}$ of $C(V\otimes_L F)$ for various ${{\beta}_1,...,{\beta}_l}$ by sign changes in front of $\sqrt{d_i}$'s and $\sqrt{-d_{m-i+1}}$'s in the definition of $f_i$ in terms of $e_j$ (see notation above). Then (if we choose all $S_{p_i}$ to be the same)
$$
(\mbox{stabilizer of}\; (p_1,...,p_r)\in {\Omega}')=(\mbox{stabilizer of}\; ({\alpha}_1,...,{\alpha}_r)\in \Omega )\cap (\mbox{stabilizer of}\; f_{1,...,1}).
$$

{\bf Remark.} Instead of $\tilde{C}(V\otimes_{L} F)$ one can also consider the Clifford algebra $C(V\otimes_{L} F)$ (or its even part ${C}^{+}(V\otimes_{L} F)$). Then the corestriction of $C(V)$ (or of $C^{+}(V)$) (with $V$ viewed as a vector space over $L$) from $L$ to $k=\mathbb Q$ (or Galois-fixed subspaces of sums (inside of tensor products of $C(V)\otimes_L F$) of tensor products of $(\Lg\otimes_k F)$-invariant $F$-vector subspaces (or ideals used above) of $C(V)\otimes_L F$, which form a single Galois orbit) would be a representation of $\Lg$ over $\mathbb Q=k$, whose extension of scalars to $F$ contains all the irreducible representations (and only them) of $\Lg\otimes_k F$ over $F$ which we need. In particular, in the case of odd $m$ it would be another primary representation of $\Lg$ over $k$.\\

\subsection{Case of the CM-field.}

Let $E=E_0(\theta), {\theta}^2\in E_0 =L$ be a CM-field.\\

Note that the tautological representation of $\Lu(\Phi)\otimes_{L,{\sigma}_i} F$ in $V\otimes_{L,{\sigma}_i} F$ splits into the direct sum of two representations of $\Lgl(m,F)\cong\Lu(\Phi)\otimes_{L,{\sigma}_i} F$:
$$
V\otimes_{L,{\sigma}_i} F=(V\otimes_{E,{\sigma}_i} F)\oplus (V\otimes_{E,\bar{{\sigma}_i}} F),
$$
where ${\sigma}_i$ and $\bar{{\sigma}_i}$ are the two extensions of ${\sigma}_i \colon E_0\rightarrow F$ to embeddings $E\rightarrow F$.\\

Since the exterior power representations ${\wedge}_F^{p}(V\otimes_{E,{\bar{{\sigma}_i}}} F)$ and ${\wedge}_F^{p}(V\otimes_{E,{{{\sigma}_i}}} F)$ of $\Lu(\Phi)\otimes_{L,{\sigma}_i} F\cong \Lgl(m,F)$ are identified by the Lie algebra automorphism $\Lgl(m,F)\rightarrow \Lgl(m, F)$, $M\mapsto -\Phi\cdot M^{T}\cdot {\Phi}^{-1}$, we have isomorphisms
$$
{\wedge}_F^{p}(V\otimes_{E,{\bar{{\sigma}_i}}} F)\rightarrow {\wedge}_F^{m-p}(V\otimes_{E,{{{\sigma}_i}}} F)\otimes_F D^{-1}
$$
and hence also isomorphisms 
$$
{\tau}_p\colon {\wedge}_F^{p}(V\otimes_{E,{\bar{{\sigma}_i}}} F)\otimes_F (E\otimes_{E,\bar{{\sigma}_i}} F)\rightarrow {\wedge}_F^{m-p}(V\otimes_{E,{{{\sigma}_i}}} F)\otimes_F D^{-1/2}, \;1\leq p\leq m
$$
of representations of $\Lgl(m,F)\cong \Lu(\Phi)\otimes_{L,{\sigma}_i} F$.\\

Let ${\wedge}_{i}^{j}\subset {\wedge}_{F}^{*}(V \otimes_{E,{\sigma}_i} F)\otimes_F F$, $1\leq i \leq r$, $1\leq j\leq m$ be the irreducible representation of $\Lgl(m,F)$ on the $F$-vector space ${\wedge}_{F}^{j}(V \otimes_{E,{\sigma}_i} F)$ twisted by $D^{-1/2}$. We define an $E_0$-linear representation $D^{c}$, $c\in\mathbb Q$ of $\Lu(\Phi)$ in the $E_0$-vector space $E$ in exactly the same way as for $\Lgl(m,F)$ above, i.e. by taking the trace of a matrix and multiplying it by $\frac{c}{m}$.\\

Consider the finite set $\Omega$ of $F$-vector spaces of the form ${\wedge}_1^{{j}_1} \otimes_F ... \otimes_F {\wedge}_r^{{j}_r}$ for various ${j}_i \in \{ 1,...,m  \}$. They are exactly the isomorphism classes of  irreducible subrepresentations of the exterior tensor product of (twisted by $D^{-1/2}$ and extended to $F$) exterior algebra representations ${\wedge}_{F}^{*}(V \otimes_{L,{\sigma}_i} F)\otimes_F (E\otimes_{L,{\sigma}_i} F)$ of $\oplus_{i=1}^{r} (\Lu(\Phi)\otimes_{L,{\sigma}_i} F)\cong \Lgl(m,F)^{\oplus r}$.\\

Let ${\wedge}_{F}^{*}(V \otimes_{L,{\sigma}_i} F)\otimes_F (E\otimes_{L,{\sigma}_i} F)=\oplus_p S_{p}^{i}$ be the decomposition into irreducible components of the representation of  $\Lu(\Phi)\otimes_{L,{\sigma}_i} F\cong \Lgl(m,F)$ obtained from the decompositions $E\otimes_{L,{\sigma}_i} F=(E\otimes_{E,{\sigma}_i} F)\oplus (E\otimes_{E,\bar{{\sigma}_i}} F)\cong D^{-1/2}\oplus D^{1/2} \cong F\oplus F$ and $V\otimes_{L,{\sigma}_i} F=(V\otimes_{E,{\sigma}_i} F)\oplus (V\otimes_{E,\bar{{\sigma}_i}} F)$ above.\\

Let ${\Omega}'$ be the finite set of $F$-vector subspaces of $({\wedge}_{F}^{*}(V \otimes_{L,{\sigma}_1} F)\otimes_F (E\otimes_{L,{\sigma}_1} F))\otimes_F ... \otimes_F ({\wedge}_{F}^{*}(V \otimes_{L,{\sigma}_r} F)\otimes_F (E\otimes_{L,{\sigma}_r} F))$ of the form ${S}_{p_1}^{1} \otimes_F ... \otimes_F {S}_{p_r}^{r}$ for various $p_1,...,p_r$. These subspaces are irreducible subrepresentations of the exterior tensor product of exterior algebra representations as a representation of $\oplus_{i=1}^{r} (\Lu(\Phi)\otimes_{L,{\sigma}_i} F)$.\\

Galois group $S=Gal(F/k)$ acts on $\Omega$ by permuting factors in tensor products. It also acts on ${\Omega}'$. Let ${\Omega}_1,...,{\Omega}_u$ be the orbits of $S$ on $\Omega$. For any $i$ choose $({j}_1,...,{j}_r)\in {\Omega}_i$ and define $U_i\subset ({\wedge}_{F}^{*}(V \otimes_{L,{\sigma}_1} F)\otimes_F (E\otimes_{L,{\sigma}_1} F))\otimes_F ... \otimes_F ({\wedge}_{F}^{*}(V \otimes_{L,{\sigma}_r} F)\otimes_F (E\otimes_{L,{\sigma}_r} F))$ to be the sum of the elements of ${\Omega}'$ (as subspaces of $({\wedge}_{F}^{*}(V \otimes_{L,{\sigma}_1} F)\otimes_F (E\otimes_{L,{\sigma}_1} F))\otimes_F ... \otimes_F ({\wedge}_{F}^{*}(V \otimes_{L,{\sigma}_r} F)\otimes_F (E\otimes_{L,{\sigma}_r} F))$) lying in the $S$-orbit of any ${S}_{p_1}^{1} \otimes_F ... \otimes_F {S}_{p_r}^{r}$, which is isomorphic to ${\wedge}_1^{{j}_1} \otimes_F ... \otimes_F {\wedge}_r^{{j}_r}$ as a $\oplus_{i=1}^{r} (\Lu(\Phi)\otimes_{L,{\sigma}_i} F)\cong {\Lgl(m,F)}^{\oplus r}$-module.\\

Then $U_i\subset ({\wedge}_{F}^{*}(V \otimes_{L,{\sigma}_1} F)\otimes_F (E\otimes_{L,{\sigma}_1} F))\otimes_F ... \otimes_F ({\wedge}_{F}^{*}(V \otimes_{L,{\sigma}_r} F)\otimes_F (E\otimes_{L,{\sigma}_r} F))$ is an $S$-submodule and
$$
\left( U_i \right)^S\cong \left( \bigoplus_{({j}_1,...,{j}_r)\in {\Omega}_i}  ({\wedge}_1^{{j}_1} \otimes_F ... \otimes_F {\wedge}_r^{{j}_r})^{\oplus n_{j_1,...,j_r}} \right)^{S}
$$
is a primary representation of $\Lg$ over $k$ of dimension $\sum_{({j}_1,...,{j}_r)\in {\Omega}_i} n_{j_1,...,j_r}\cdot {\binom{m}{j_1}} \cdot ... \cdot {\binom{m}{j_r}}$. These representations $\left( U_i \right)^S$, $1\leq i \leq u$ contain all representations of $\Lg\otimes_k F$ of the form ${\wedge}_1^{{j}_1} \otimes_F ... \otimes_F {\wedge}_r^{{j}_r}$ after extending scalars to $F$.\\

The reason why nontrivial multiplicities may appear is exactly the doubling $V\otimes_{L,{\sigma}_i} F=(V\otimes_{E,{\sigma}_i} F)\oplus (V\otimes_{E,\bar{{\sigma}_i}} F)$ described above. Hence one can compute multiplicities $n_{j_1,...,j_r}$ as follows. Consider the finite set ${\Omega}''$ of $r$-tuples of signs $+$ and $-$, i.e. ${\Omega}''=\{  ({\alpha}_1,...,{\alpha}_r)\;\mid\; {\alpha}_i=\pm \}$. Note that the $i$-th sign corresponds to the $i$-th embedding ${\sigma}_i\colon L\rightarrow F$ over $k$. Consider the action of $S=Gal(F/k)$ on ${\Omega}''$ such that $g\in S$ acts on entries of $r$-tuples by the same permutations as on the set of left cosets $S/\tilde{H}$ (where $\tilde{H}=\{ g\in S \;\mid \; g\circ {\sigma}_1={\sigma}_1 \}$) and $g$ changes the sign in the $i$-th entry to the opposit sign (in the $j$-th entry, where ${\sigma}_j=g\circ {\sigma}_i$) if and only if $g(\theta)=-\theta$. Then
$$
n_{j_1,...,j_r}=\frac{\mbox{order of the stabilizer of}\; ({j}_1,...,{j}_r)\in\Omega}{\mbox{order of the intersection of stabilizers of}\; (+,...,+)\in{\Omega}''\;\mbox{ and of }\; ({j}_1,...,{j}_r)\in \Omega }.
$$

This gives a description of some multiples of $(k=\mathbb Q)$-linear irreducible representations $W_i$ of $\Lg$ mentioned in the Theorem above (as well as formulas for their dimensions - some multiples of $dim_k(W_i)$) in terms of the Galois action.\\

\section{Cohomology classes of division algebras.} 

In this section we compute division algebras $D_i$ as elements of the Brauer group $Br(F/C_j)\cong H^2(Gal(F/C_j), F^{*})$ as well as their centers $C_j$.\\ 

\subsection{Case of the totally real field and odd dimension.}

Let $E=E_0=L$ be totally real and $m=dim_EV$ odd. We saw above how to construct a primary representation $W=U^S$ of $\Lg$ over $k=\mathbb Q$, which contains irreducible representation ${\rho}^0 \boxtimes ... \boxtimes {\rho}^0$ (the exterior tensor product of irreducible spin representations) of $\Lg\otimes_k F\cong \oplus_{i=1}^{r} \Lso(\Phi)\otimes_{L,{\sigma}_i} F$ after extending scalars to $F$. This means that $W\cong W_0^{\oplus \mu}$, where $W_0$ is an irreducible representation of $\Lg$ over $k$ and $W_0\otimes_k F\cong \frac{dim_kW}{\mu \cdot (dim_F({\rho}^0))^r}\cdot {\rho}^0 \boxtimes ... \boxtimes {\rho}^0$. Since we are interested only in the endomorphism algebra $D_0=End_{\Lg}(W_0)$ which is a central division algebra over $k$ split over $F$, we can describe it by computing the Galois cohomology invariant of the central simple algebra $A=End_{\Lg}(W)\cong Mat_{\mu \times \mu}(D_0)$, i.e. its Brauer invariant in $Br(F/k)\cong H^2(S,F^{*})$, where $S=Gal(F/k)$. Then $\mu=\frac{deg(A)}{deg(D_0)}=\frac{n_0}{deg(D_0)}$.\\

We will use the same notation as above with the following exceptions:
$$
f_{{\alpha}_1,...,{\alpha}_l,\gamma} = (1+\gamma \cdot f_0)\cdot f_{{\alpha}_1\cdot 1}\cdot ... \cdot f_{{\alpha}_l\cdot l},
$$

$$
f_{{\alpha}\cdot i}=\left(e_i+\alpha \cdot \frac{\sqrt{d_i}}{\sqrt{-d_{m-i+1}}}\cdot e_{m-i+1}\right).
$$
Some parts of our construction (in particular, the construction of the generators of endomorphism algebras) may be viewed as a generalization of some constructions of van Geemen \cite{vanGeemen1}, \S 3.\\

Consider $F$-linear homomorphisms 
$$
r_{(({\alpha}_i),\gamma),(({\beta}_i),\tilde{\gamma})} \colon \tilde{C}(V\otimes_L F)\rightarrow I_{{\beta}_1,...,{\beta}_l,\tilde{\gamma}}, \;\xi \mapsto {\tau}^{{\delta}(\gamma,\tilde{\gamma})} (\xi\cdot R_{(({\alpha}_i),\gamma),({\beta}_i)}),
$$
where $\tau \colon {C}(V\otimes_L F)\rightarrow {C}(V\otimes_L F)$ is the algebra homomorphism induced by multiplication by $(-1)$ on $V$, ${\delta}(\gamma,\tilde{\gamma}) = 1$, if $\gamma\neq \tilde{\gamma}\cdot (-1)^{P(\alpha,\beta)}$ (where $P(\alpha,\beta)=card \{ i \;\mid\; {\alpha}_i\neq {\beta}_i \}$) and $0$ otherwise, and
$$
R_{(({\alpha}_i),\gamma),({\beta}_i)}=\frac{(-1)^{c(\alpha,\beta)}}{\prod_{i\colon {\alpha}_i= {\beta}_i} \Phi (f_i,f_{-i})} \cdot \prod_{i\colon {\alpha}_i= {\beta}_i} (f_{-{\alpha}_i \cdot i}\cdot f_{{\alpha}_i \cdot i}) \cdot \prod_{i\colon {\alpha}_i\neq {\beta}_i} f_{{\beta}_i \cdot i},
$$
where $c(\alpha,\beta)$ is the number of transpositions of factors needed to transform the product $\prod_i f_{{\alpha}_i \cdot i}\cdot \prod_{i\colon {\alpha}_i\neq {\beta}_i} f_{{\beta}_i \cdot i}$ into the product $q\cdot \prod_i f_{{\beta}_i \cdot i}$ with some coefficient $q\in C(V\otimes_L F)$. Then $r_{(({\alpha}_i),\gamma),(({\beta}_i),\tilde{\gamma})}$ is nonzero only on the factor $I_{{\alpha}_1,...,{\alpha}_l,{\gamma}}$ of $\tilde{C}(V\otimes_L F)$ and induces an isomorphism $I_{{\alpha}_1,...,{\alpha}_l,{\gamma}} \rightarrow I_{{\beta}_1,...,{\beta}_l,\tilde{\gamma}}$ which commutes with the action of $\Lso(\Phi)\otimes_L F$.\\

In order to simplify notation we will denote index $(({\alpha}_i),\gamma)$ by $\alpha$.\\

One can choose coefficients ${\lambda}_{\alpha,\beta}\in F^{*}$ such that under an isomorphism of $F$-algebras $Mat(F)\cong End_{\Lso(\Phi)\otimes_L F}(\tilde{C}(V\otimes_L F))$ matrices of the form $E_{ij}$ (in the notation of \cite{Bourbaki}, \S 13) correspond to endomorphisms ${\lambda}_{\alpha,\beta}\cdot r_{\alpha,\beta}$. In order to do this, one can choose and fix index ${\alpha}^0=(({\alpha}_i^0),{\gamma}^0)$ and take
$$
{\lambda}_{{\alpha}^0,\beta}=1, \; {\lambda}_{\beta,{\alpha}^0}=(-1)^{P({\alpha}^0,\beta)\cdot {\delta}({\gamma}^0,\tilde{\gamma})+P({\alpha}^0,\beta)\cdot (P({\alpha}^0,\beta)-1)/2}\cdot \prod_{i\colon {\alpha}_i^0 \neq {\beta}_i}\frac{1}{\Phi(f_i,f_{-i})}
$$
and
$$
{\lambda}_{\alpha,\beta}={\lambda}_{\alpha,{\alpha}^0}\cdot (-1)^{e(\alpha,\beta)+ {\delta}(\gamma,\tilde{\gamma})\cdot (l+P(\alpha,\beta))   + {\delta}(\gamma,{\gamma}^0)\cdot (l+P(\alpha,{\alpha}^0))   + {\delta}({\gamma}^0,\tilde{\gamma})\cdot (l+P({\alpha}^0,\beta))  }\cdot \prod_{i\colon {\alpha}_i = {\beta}_i \neq {\alpha}_i^0  } \Phi(f_i,f_{-i}),
$$
where ${\alpha}=(({\alpha}_i),{\gamma})$, ${\beta}=(({\beta}_i),\tilde{\gamma})$, $e(\alpha,\beta)$ is the number of transpositions of factors needed in order to transform the product $\prod_{i\colon {\alpha}_i^0 \neq {\alpha}_i} f_{{\alpha}_i\cdot i}\cdot \prod_{i\colon {\alpha}_i^0 \neq {\beta}_i} f_{-{\beta}_i\cdot i}$ into the product $\prod_{i\colon {\alpha}_i \neq {\beta}_i} f_{{\alpha}_i\cdot i}\cdot \prod_{i\colon {\alpha}_i= {\beta}_i \neq {\alpha}_i^0} (f_{{\beta}_i\cdot i} \cdot f_{-{\beta}_i\cdot i} )$. Note that in this construction ${\lambda}_{\alpha,\beta}\in L^{*}$.\\

Then we construct endomorphisms 
\begin{multline*}
r_{({\alpha}^i),({\beta}^i)}=r_{{\alpha}^1,{\beta}^1}^{1}\circ ... \circ r_{{\alpha}^r,{\beta}^r}^{r}\colon \tilde{C}(V\otimes_{L,{\sigma}_1} F)\otimes_F ... \otimes_F \tilde{C}(V\otimes_{L,{\sigma}_r} F)\rightarrow \\
\rightarrow I_{{\beta}^1}^1\otimes_F ... \otimes_F I_{{\beta}^r}^r\subset \tilde{C}(V\otimes_{L,{\sigma}_1} F)\otimes_F ... \otimes_F \tilde{C}(V\otimes_{L,{\sigma}_r} F)
\end{multline*}
which commute with $\Lg\otimes_k F$, where ${\alpha}^p=(({\alpha}_{1}^{p},...,{\alpha}_{l}^{p}),{\gamma}^p)$, ${\beta}^p=(({\beta}_{1}^{p},...,{\beta}_{l}^{p}),{\tilde{\gamma}}^p)$ and
\begin{multline*}
r_{{\alpha}^p,{\beta}^p}^{p}=1\otimes_F ... \otimes_F (r_{{\alpha}^p,{\beta}^p}) \otimes_F ... \otimes_F 1\colon \tilde{C}(V\otimes_{L,{\sigma}_1} F)\otimes_F ... \otimes_F \tilde{C}(V\otimes_{L,{\sigma}_r} F)\rightarrow \\
\rightarrow \tilde{C}(V\otimes_{L,{\sigma}_1} F)\otimes_F ... \otimes_F \tilde{C}(V\otimes_{L,{\sigma}_r} F)
\end{multline*}
(with $1$ outside of the $p$-th spot).\\

As in \cite{vanGeemen1}, Proposition 3.6 $F$-algebra $End_{\Lg\otimes_k F}(W\otimes_k F)=A\otimes_k F$ is generated by elements $r_{({\alpha}^i),({\beta}^i)}$ (more precisely, by those of them which correspond to the summands of $\tilde{C}(V\otimes_{L,{\sigma}_1} F)\otimes_F ... \otimes_F \tilde{C}(V\otimes_{L,{\sigma}_r} F)$ included in $W\otimes_k F=U\subset \tilde{C}(V\otimes_{L,{\sigma}_1} F)\otimes_F ... \otimes_F \tilde{C}(V\otimes_{L,{\sigma}_r} F)$) or by elements $r_{{\alpha},{\beta}}^{p}$, while $k$-algebra $A=End_{\Lg}(W)=(A\otimes_k F)^S$ is generated by elements $r_{{\alpha},{\beta}}^{p,q}=\sum_{g \in S} g(e_q)\cdot g\circ r_{{\alpha},{\beta}}^{p}$, where $\{  e_q \}$ is a basis of $F/k$.\\

Let us denote by $(c_{q,g})$ the inverse matrix of the matrix $(g(e_q))$. Then $r_{\alpha,\beta}^p=\sum_{q} c_{q,Id}\cdot r_{\alpha,\beta}^{p,q}$ and for any $g\in S=Gal(F/k)$ if we denote by ${\phi}_g\colon A\otimes_k F\rightarrow A\otimes_k F$ the conjugation by $g\colon a\otimes f\mapsto a\otimes g(f)$, then
$$
{\phi}_g(r_{({\alpha}^i),({\beta}^i)})=g\circ r_{({\alpha}^i),({\beta}^i)}= r_{g({\alpha}^i),g({\beta}^i)},
$$
where the action of $S$ on upper indices $i$ (which number embeddings ${\sigma}_i\colon L\hookrightarrow F$) coincides with its action on left cosets $S/ \tilde{H}$, where $\tilde{H}=\{ g\in S \; \mid\; g {\mid}_{{\sigma}_1(L)}=Id_{{\sigma}_1(L)} \}$ and the action of $g\in S$ on indices $\alpha=(({\alpha}_1,...,{\alpha}_l),\gamma)$ is given by the rule $g(\alpha)=((c_1(g)\cdot{\alpha}_1,...,c_l(g)\cdot{\alpha}_l),c_0(g)\cdot \gamma)$, where $c_i(g)\in \{ \pm 1 \}$ and $g(f_{{\alpha}_1,...,{\alpha}_l,\gamma})=f_{c_1(g)\cdot{\alpha}_1,...,c_l(g)\cdot{\alpha}_l,c_0(g)\cdot \gamma}$.\\

Hence the matrix of $m(g)\in GL(W\otimes_k F)$ is such that
$$
m(g)\cdot E_{i,j}\cdot {m(g)}^{-1}={\phi}_g (E_{i,j})=\left(  \prod_{i=1}^r \frac{g({\lambda}_{{\alpha}^i,{\beta}^i})}{{\lambda}_{g({\alpha}^i),g({\beta}^i)}} \right)  \cdot E_{g(i),g(j)},
$$
where $E_{i,j}$ denotes a matrix from $Mat(F)\cong End_{\Lg\otimes_k F} (W\otimes_k F)$ corresponding to $r_{({\alpha}^i),({\beta}^i)}$, i.e. upto a scalar multiple conjugation by $m(g)$ acts on matrices as the (same) permutation of columns and rows induced by $g$ on indices $(({\alpha}_{j}^{i}),{\gamma}^i)$.\\

Then the element of $H^2(S,F^{*})$ corresponding to the central division algebra $D_0=End_{\Lg}(W_0)$ is the class of a $2$-cocycle $\lambda\colon S\times S\rightarrow F^{*}\cong F^{*}\cdot Id\subset Mat(F)$, $(g_1,g_2)\mapsto m(g_1g_2)\cdot (g_1 (m(g_2)))^{-1} \cdot m(g_1)^{-1}$ \cite{Kuznetsov}, \cite{Jacobson}.\\

\subsection{Case of the totally real field and even dimension.}

Let $E=E_0=L$ be totally real and $m=dim_EV$ even. We saw above how to construct a primary representation $W=(U_i)^S$ of $\Lg$ over $k=\mathbb Q$, which contains irreducible representation ${\rho}^{{\alpha}_1} \boxtimes ... \boxtimes {\rho}^{{\alpha}_r}$ (the exterior tensor product of irreducible semi-spin representations) of $\Lg\otimes_k F\cong \oplus_{i=1}^{r} \Lso(\Phi)\otimes_{L,{\sigma}_i} F$ after extending scalars to $F$ (as well as its Galois conjugates). This means that $W\cong W_0^{\oplus \mu}$, where $W_0$ is an irreducible representation of $\Lg$ over $k$, $W\otimes_k F\cong \oplus_i W_i$ and $W_i\cong \frac{dim_F{W_i}}{(dim_F({\rho}^{{\alpha}_1}))^r}\cdot {\rho}^{{{\alpha}_1}'} \boxtimes ... \boxtimes {\rho}^{{{\alpha}_r}'}$ are the isotypical components (over $F$). Since we are interested only in the endomorphism algebra $D_0=End_{\Lg}(W_0)$ which is a division algebra over $k$ (and over its center $C$) split over $F$, we can describe it by computing the Galois cohomology invariant of the central simple algebra $A=End_{\Lg}(W)\cong Mat_{\mu \times \mu}(D_0)$ (over $C$), i.e. its Brauer invariant in $Br(F/C)\cong H^2(S',F^{*})$, where $S'=Gal(F/C)$. Then $\mu=\frac{deg(A)}{deg(D_0)}=\frac{n_{{\alpha}_1,...,{\alpha}_r}}{deg(D_0)}$.\\

We will use the same notation as above with the following exceptions:
$$
f_{{\alpha}_1,...,{\alpha}_l} = f_{{\alpha}_1\cdot 1}\cdot ... \cdot f_{{\alpha}_l\cdot l},
$$

$$f_{{\alpha}\cdot i}=\left(e_i+\alpha \cdot \frac{\sqrt{d_i}}{\sqrt{-d_{m-i+1}}}\cdot e_{m-i+1}\right).
$$
Some parts of our construction (in particular, the construction of the generators of endomorphism algebras) may be viewed as a generalization of some constructions of van Geemen \cite{vanGeemen1}, \S 3.\\

Consider $F$-linear homomorphisms 
$$
r_{({\alpha}_i),({\beta}_i)} \colon \tilde{C}(V\otimes_L F)\rightarrow I_{{\beta}_1,...,{\beta}_l}, \; \xi \mapsto \xi\cdot R_{({\alpha}_i),({\beta}_i)},
$$
where $P(\alpha,\beta)=card \{ i \;\mid\; {\alpha}_i\neq {\beta}_i \}$ and 
$$
R_{({\alpha}_i),({\beta}_i)}=\frac{(-1)^{c(\alpha,\beta)}}{\prod_{i\colon {\alpha}_i= {\beta}_i} \Phi (f_i,f_{-i})} \cdot \prod_{i\colon {\alpha}_i= {\beta}_i} (f_{-{\alpha}_i \cdot i}\cdot f_{{\alpha}_i \cdot i}) \cdot \prod_{i\colon {\alpha}_i\neq {\beta}_i} f_{{\beta}_i \cdot i},
$$
where $c(\alpha,\beta)$ is the number of transpositions of factors needed to transform the product $\prod_i f_{{\alpha}_i \cdot i}\cdot \prod_{i\colon {\alpha}_i\neq {\beta}_i} f_{{\beta}_i \cdot i}$ into the product $q\cdot \prod_i f_{{\beta}_i \cdot i}$ with some coefficient $q\in C(V\otimes_L F)$. Then $r_{({\alpha}_i),({\beta}_i)}$ is nonzero only on the factor $I_{{\alpha}_1,...,{\alpha}_l}$ of $\tilde{C}(V\otimes_L F)$ and induces an isomorphism $I_{{\alpha}_1,...,{\alpha}_l} \rightarrow I_{{\beta}_1,...,{\beta}_l}$ which commutes with action of $\Lso(\Phi)\otimes_L F$. Without mentioning this explicitely, we will be restricting all our endomorphisms to the factors of $\tilde{C}(V\otimes_{L,{\sigma}_1} F)\otimes_F ... \otimes_F \tilde{C}(V\otimes_{L,{\sigma}_r} F)$ contributing to an isotypical component $W_i\subset \tilde{C}(V\otimes_{L,{\sigma}_1} F)\otimes_F ... \otimes_F \tilde{C}(V\otimes_{L,{\sigma}_r} F)$.\\

In order to simplify notation we will denote index $({\alpha}_i)$ by $\alpha$.\\

One can choose coefficients ${\lambda}_{\alpha,\beta}\in F^{*}$ such that under an isomorphism of $F$-algebras $Mat(F)\cong End_{\Lso(\Phi)\otimes_L F}(W_i)$ (note that $W_i\subset \tilde{C}(V\otimes_{L,{\sigma}_1} F)\otimes_F ... \otimes_F \tilde{C}(V\otimes_{L,{\sigma}_r} F)$ and see the remark above) matrices of the form $E_{ij}$ correspond to endomorphisms ${\lambda}_{\alpha,\beta}\cdot r_{\alpha,\beta}$. In order to do this, one can choose and fix index ${\alpha}^0=({\alpha}_i^0)$ and take
$$
{\lambda}_{{\alpha}^0,\beta}=1, \; {\lambda}_{\beta,{\alpha}^0}=(-1)^{P({\alpha}^0,\beta)\cdot (P({\alpha}^0,\beta)-1)/2}\cdot \prod_{i\colon {\alpha}_i^0 \neq {\beta}_i}\frac{1}{\Phi(f_i,f_{-i})}
$$
and
$$
{\lambda}_{\alpha,\beta}={\lambda}_{\alpha,{\alpha}^0}\cdot (-1)^{e(\alpha,\beta)}\cdot \prod_{i\colon {\alpha}_i = {\beta}_i \neq {\alpha}_i^0  } \Phi(f_i,f_{-i}),
$$
where ${\alpha}=({\alpha}_i)$, ${\beta}=({\beta}_i)$, $e(\alpha,\beta)$ is the number of transpositions of factors needed in order to transform the product $\prod_{i\colon {\alpha}_i^0 \neq {\alpha}_i} f_{{\alpha}_i\cdot i}\cdot \prod_{i\colon {\alpha}_i^0 \neq {\beta}_i} f_{-{\beta}_i\cdot i}$ into the product $\prod_{i\colon {\alpha}_i \neq {\beta}_i} f_{{\alpha}_i\cdot i}\cdot \prod_{i\colon {\alpha}_i= {\beta}_i \neq {\alpha}_i^0} (f_{{\beta}_i\cdot i} \cdot f_{-{\beta}_i\cdot i} )$. Note that in this construction ${\lambda}_{\alpha,\beta}\in L^{*}$.\\

Then we construct endomorphisms 
\begin{multline*}
r_{({\alpha}^i),({\beta}^i)}=r_{{\alpha}^1,{\beta}^1}^{1}\circ ... \circ r_{{\alpha}^r,{\beta}^r}^{r}\colon \tilde{C}(V\otimes_{L,{\sigma}_1} F)\otimes_F ... \otimes_F \tilde{C}(V\otimes_{L,{\sigma}_r} F)\rightarrow \\
\rightarrow I_{{\beta}^1}^1\otimes_F ... \otimes_F I_{{\beta}^r}^r\subset \tilde{C}(V\otimes_{L,{\sigma}_1} F)\otimes_F ... \otimes_F \tilde{C}(V\otimes_{L,{\sigma}_r} F)
\end{multline*}
which commute with $\Lg\otimes_k F$, where ${\alpha}^p=({\alpha}_{1}^{p},...,{\alpha}_{l}^{p})$, ${\beta}^p=({\beta}_{1}^{p},...,{\beta}_{l}^{p})$ and
\begin{multline*}
r_{{\alpha}^p,{\beta}^p}^{p}=1\otimes_F ... \otimes_F (r_{{\alpha}^p,{\beta}^p}) \otimes_F ... \otimes_F 1\colon \tilde{C}(V\otimes_{L,{\sigma}_1} F)\otimes_F ... \otimes_F \tilde{C}(V\otimes_{L,{\sigma}_r} F)\rightarrow \\
\rightarrow \tilde{C}(V\otimes_{L,{\sigma}_1} F)\otimes_F ... \otimes_F \tilde{C}(V\otimes_{L,{\sigma}_r} F)
\end{multline*}
(with $1$ outside of the $p$-th spot).\\

As in \cite{vanGeemen1}, Proposition 3.6 $F$-algebra $End_{\Lg\otimes_k F}(W\otimes_k F)=A\otimes_k F$ is generated by elements $r_{({\alpha}^i),({\beta}^i)}$ (more precisely, by those of them which correspond to the summands of $\tilde{C}(V\otimes_{L,{\sigma}_1} F)\otimes_F ... \otimes_F \tilde{C}(V\otimes_{L,{\sigma}_r} F)$ included in various isotypical components $W_{i'}\otimes_k F\subset U_i \subset \tilde{C}(V\otimes_{L,{\sigma}_1} F)\otimes_F ... \otimes_F \tilde{C}(V\otimes_{L,{\sigma}_r} F)$) or by elements $r_{{\alpha},{\beta}}^{p}$, while $k$-algebra $A=End_{\Lg}(W)=(A\otimes_k F)^S$ is generated by elements $r_{{\alpha},{\beta}}^{p,q}=\sum_{g \in S} g(e_q)\cdot g\circ r_{{\alpha},{\beta}}^{p}$, where $\{  e_q \}$ is a basis of $F/k$.\\

The center $C$ of $A$ (and of $D_0$) consists of Galois averages (as above) of $F$-linear combinations of sums $C_{{i}'}=\sum_{({\alpha}^{j})\in I_{{i}'}} (\prod_{i=1}^{r} {{\sigma}_i}({\lambda}_{{\alpha}^i,{\alpha}^i})) \cdot r_{({\alpha}^j),({\alpha}^j)}$ (over the sets $I_{{i}'}$ of indices ${\alpha}^{j}$ corresponding to irreducible subrepresentations over $F$ of $W\otimes_k F$ contained in various isotypical components $W_{i'}$). Each of the coefficients of these $F$-linear combinations gives a field embedding $C\rightarrow F$ over $k=\mathbb Q$. Note that $A\otimes_k F \cong \prod A\otimes_C F$, where the product is taken over these embeddings (which are numbered by the isotypical components $W_{{i}'}$ of $W\otimes_k F$ over $F$) and $A\otimes_C F\cong End_{\Lg\otimes_k F}(W_{{i}'})$. Moreover, the projection $A\otimes_k F \rightarrow A\otimes_C F$ is given by annihilating endomorphisms between irreducible subrepresentations of isotypical components $W_{{i}''}$ different from $W_{{i}'}$. More explicitely the subfield $C\subset F$ under the embedding corresponding to an isotypical component $W_{{i}'}$ is the fixed subfield of the subgroup $S'\subset S$ consisting of those $g\in S$ which preserve the isotypical component: $g(W_{{i}'})=W_{{i}'}$. Let us choose one such embedding $C\rightarrow F$ (which corresponds to a choice of an isotypical component $W_{{i}'}$ of $W\otimes_k F$).\\

Let us denote by $(c_{q,g})$ the inverse matrix of the matrix $(g(e_q))$. Then $r_{\alpha,\beta}^p=\sum_{q} c_{q,Id}\cdot r_{\alpha,\beta}^{p,q}$ and for any $g\in S'=Gal(F/C)\subset S=Gal(F/k)$ if we denote by ${\phi}_g\colon A\otimes_C F\rightarrow A\otimes_C F$ the conjugation by $g\colon a\otimes f\mapsto a\otimes g(f)$, then
$$
{\phi}_g(r_{({\alpha}^i),({\beta}^i)})=g\circ r_{({\alpha}^i),({\beta}^i)}= r_{g({\alpha}^i),g({\beta}^i)},
$$
where the action of $S'\subset S$ on upper indices $i$ (which number embeddings ${\sigma}_i\colon L\hookrightarrow F$) coincides with its action on left cosets $S/ \tilde{H}$, where $\tilde{H}=\{ g\in S \; \mid\; g {\mid}_{{\sigma}_1(L)}=Id_{{\sigma}_1(L)} \}$ and the action of $g\in S'\subset S$ on indices $\alpha=({\alpha}_1,...,{\alpha}_l)$ is given by the rule $g(\alpha)=(c_1(g)\cdot{\alpha}_1,...,c_l(g)\cdot{\alpha}_l)$, where $c_i(g)\in \{ \pm 1 \}$ and $g(f_{{\alpha}_1,...,{\alpha}_l})=f_{c_1(g)\cdot{\alpha}_1,...,c_l(g)\cdot{\alpha}_l}$.\\

Hence the matrix of $m(g)\in GL(W_{{i}'})$ is such that 
$$
m(g)\cdot E_{i,j}\cdot {m(g)}^{-1}={\phi}_g (E_{i,j})=\left(  \prod_{i=1}^r \frac{g({\lambda}_{{\alpha}^i,{\beta}^i})}{{\lambda}_{g({\alpha}^i),g({\beta}^i)}} \right)  \cdot E_{g(i),g(j)},
$$
where $E_{i,j}$ denotes a matrix from $Mat(F)\cong End_{\Lg\otimes_k F} (W_{{i}'})$ corresponding to $r_{({\alpha}^i),({\beta}^i)}$, i.e. upto a scalar multiple conjugation by $m(g)$ acts on matrices as the (same) permutation of columns and rows induced by $g$ on indices $({\alpha}_{j}^{i})$.\\

Then the element of $H^2(S',F^{*})$ corresponding to the central division algebra $D_0=End_{\Lg}(W_0)$ (over $C$) is the class of a $2$-cocycle $\lambda\colon S'\times S'\rightarrow F^{*}\cong F^{*}\cdot Id\subset Mat(F)$, $(g_1,g_2)\mapsto m(g_1g_2)\cdot (g_1 (m(g_2)))^{-1} \cdot m(g_1)^{-1}$ \cite{Kuznetsov}, \cite{Jacobson}.\\

\subsection{Case of the CM-field.}

Let $E=E_0(\theta), {\theta}^2\in E_0=L$ be a CM-field and $m=dim_EV$. We saw above how to construct a primary representation $W=(U_i)^S$ of $\Lg$ over $k=\mathbb Q$, which contains the irreducible representation ${\rho}_{j_1}^{{\alpha}_1} \boxtimes ... \boxtimes {\rho}_{j_r}^{{\alpha}_r}$ of $\Lg\otimes_k F\cong \oplus_{i=1}^{r} \Lu(\Phi)\otimes_{L,{\sigma}_i} F\cong {\Lgl(m,F)}^{\oplus r}$ after extending scalars to $F$ (as well as its Galois conjugates). Here ${\alpha}_{i}\in \{ \pm \}$, $1\leq j_i\leq m$ and
$$
{\rho}_{j_i}^{{\alpha}_i}\colon \Lgl(m,F)\rightarrow End_{F}({\wedge}_{F}^{j_i}(V\otimes_{E,{\alpha_i\cdot {\sigma}}} F)\otimes_F F)
$$
is the exterior product representation twisted by $D^{{\alpha}_i/2}$, where $\pm {\sigma}\colon E\rightarrow F$ are the two embeddings extending ${\sigma}\colon L\rightarrow F$. This means that $W\cong W_0^{\oplus \mu}$, where $W_0$ is an irreducible representation of $\Lg$ over $k$, $W\otimes_k F\cong \oplus_i W_i$ and $W_i\cong \frac{dim_F{W_i}}{dim_F({\rho}_{j_1}^{{\alpha}_1})\cdot ... \cdot dim_F({\rho}_{j_r}^{{\alpha}_r})}\cdot  {\rho}_{j_1}^{{{\alpha}_1}'} \boxtimes ... \boxtimes {\rho}_{j_r}^{{{\alpha}_r}'}$ are the isotypical components (over $F$). Since we are interested only in the endomorphism algebra $D_0=End_{\Lg}(W_0)$ which is a division algebra over $k$ (and over its center $C$) split over $F$, we can describe it by computing the Galois cohomology invariant of the central simple algebra $A=End_{\Lg}(W)\cong Mat_{\mu \times \mu}(D_0)$ (over $C$), i.e. its Brauer invariant in $Br(F/C)\cong H^2(S',F^{*})$, where $S'=Gal(F/C)$. Then $\mu=\frac{deg(A)}{deg(D_0)}=\frac{n_{{j}_1,...,{j}_r}}{deg(D_0)}$.\\

Our computation is analogous to the case of a totally real field considered above.\\

Consider $F$-linear homomorphisms 
$$
r_{{\alpha},{\beta}} \colon {\wedge}_{F}^{*}(V\otimes_{E,{\alpha\cdot {\sigma}}} F)\otimes_F F\rightarrow {\wedge}_{F}^{*}(V\otimes_{E,{\beta\cdot {\sigma}}} F)\otimes_F F, \; \xi\mapsto ({\tau}_{*})^{P(\alpha,\beta)} (\xi),
$$
where $P(-1,+1)=1$, $P(+1,-1)=-1$, $P(\alpha,\alpha)=0$ and ${\tau}_{*}=\oplus_p {\tau}_p$ is the direct sum of isomorphisms of $\Lgl(m,F)$-modules
$$
{\wedge}_F^{p}(V\otimes_{E,{\bar{{\sigma}_i}}} F)\otimes_F (E\otimes_{E,\bar{{\sigma}_i}} F)\rightarrow {\wedge}_F^{m-p}(V\otimes_{E,{{{\sigma}_i}}} F)\otimes_F D^{-1/2}
$$
introduced above. Then $r_{{\alpha},{\beta}}$ induces an isomorphism 
$$
{\wedge}_{F}^{j_i}(V\otimes_{E,{\alpha\cdot {\sigma}}} F)\otimes_F F\rightarrow {\wedge}_{F}^{{j_i}'}(V\otimes_{E,{\beta\cdot {\sigma}}} F)\otimes_F F
$$
which commutes with the action of $\Lu(\Phi)\otimes_L F\cong \Lgl(m,F)$. Without mentioning this explicitely, we will be restricting all our endomorphisms to the factors of $({\wedge}_{F}^{*}(V \otimes_{L,{\sigma}_1} F)\otimes_F (E\otimes_{L,{\sigma}_1} F))\otimes_F ... \otimes_F ({\wedge}_{F}^{*}(V \otimes_{L,{\sigma}_r} F)\otimes_F (E\otimes_{L,{\sigma}_r} F))$ contributing to an isotypical component $W_i\subset ({\wedge}_{F}^{*}(V \otimes_{L,{\sigma}_1} F)\otimes_F (E\otimes_{L,{\sigma}_1} F))\otimes_F ... \otimes_F ({\wedge}_{F}^{*}(V \otimes_{L,{\sigma}_r} F)\otimes_F (E\otimes_{L,{\sigma}_r} F))$.\\

Then we construct endomorphisms 
\begin{multline*}
r_{({\alpha}^i),({\beta}^i)}=r_{{\alpha}^1,{\beta}^1}^{1}\circ ... \circ r_{{\alpha}^r,{\beta}^r}^{r}\colon ({\wedge}_{F}^{*}(V \otimes_{L,{\sigma}_1} F)\otimes_F (E\otimes_{L,{\sigma}_1} F))\otimes_F \ldots \\
\ldots \otimes_F ({\wedge}_{F}^{*}(V \otimes_{L,{\sigma}_r} F)\otimes_F (E\otimes_{L,{\sigma}_r} F))\rightarrow \\
\rightarrow ({\wedge}_{F}^{*}(V \otimes_{L,{\sigma}_1} F)\otimes_F (E\otimes_{L,{\sigma}_1} F))\otimes_F ... \otimes_F ({\wedge}_{F}^{*}(V \otimes_{L,{\sigma}_r} F)\otimes_F (E\otimes_{L,{\sigma}_r} F)),
\end{multline*}
which commute with $\Lg\otimes_k F$, where $({\alpha}^i)=({\alpha}^{1},...,{\alpha}^{r})$, $({\beta}^i)=({\beta}^{1},...,{\beta}^{r})$ and
\begin{multline*}
r_{{\alpha}^p,{\beta}^p}^{p}=1\otimes_F ... \otimes_F (r_{{\alpha}^p,{\beta}^p})\otimes_F ... \otimes_F 1\colon ({\wedge}_{F}^{*}(V \otimes_{L,{\sigma}_1} F)\otimes_F (E\otimes_{L,{\sigma}_1} F))\otimes_F \ldots \\
\ldots \otimes_F ({\wedge}_{F}^{*}(V \otimes_{L,{\sigma}_r} F)\otimes_F (E\otimes_{L,{\sigma}_r} F))\rightarrow \\
\rightarrow ({\wedge}_{F}^{*}(V \otimes_{L,{\sigma}_1} F)\otimes_F (E\otimes_{L,{\sigma}_1} F))\otimes_F ... \otimes_F ({\wedge}_{F}^{*}(V \otimes_{L,{\sigma}_r} F)\otimes_F (E\otimes_{L,{\sigma}_r} F))
\end{multline*}
(with $1$ outside of the $p$-th spot).\\

As in the case of a totally real field $E$, $F$-algebra $End_{\Lg\otimes_k F}(W\otimes_k F)=A\otimes_k F$ is generated by elements $r_{({\alpha}^i),({\beta}^i)}$ (more precisely, by those of them which correspond to the summands of $({\wedge}_{F}^{*}(V \otimes_{L,{\sigma}_1} F)\otimes_F (E\otimes_{L,{\sigma}_1} F))\otimes_F ... \otimes_F ({\wedge}_{F}^{*}(V \otimes_{L,{\sigma}_r} F)\otimes_F (E\otimes_{L,{\sigma}_r} F))$ included in various isotypical components $W_{i'}\otimes_k F\subset U_i \subset ({\wedge}_{F}^{*}(V \otimes_{L,{\sigma}_1} F)\otimes_F (E\otimes_{L,{\sigma}_1} F))\otimes_F ... \otimes_F ({\wedge}_{F}^{*}(V \otimes_{L,{\sigma}_r} F)\otimes_F (E\otimes_{L,{\sigma}_r} F))$) or by elements $r_{{\alpha},{\beta}}^{p}$, while $k$-algebra $A=End_{\Lg}(W)=(A\otimes_k F)^S$ is generated by elements $r_{{\alpha},{\beta}}^{p,q}=\sum_{g \in S} g(e_q)\cdot g\circ r_{{\alpha},{\beta}}^{p}$, where $\{  e_q \}$ is a basis of $F/k$.\\

The center $C$ of $A$ (and of $D_0$) can be computed exactly as in the case of a totally real field. In particular, field embeddings $C\rightarrow F$ correspond to the isotypical components $W_{{i}'}$ of $W\otimes_k F$ over $F$, $A\otimes_C F\cong End_{\Lg\otimes_k F}(W_{{i}'})$, the projection $A\otimes_k F\cong \prod A\otimes_C F \rightarrow A\otimes_C F$ is given by annihilating endomorphisms between irreducible subrepresentations of isotypical components $W_{{i}''}$ different from $W_{{i}'}$ and the subfield $C\subset F$ under the embedding corresponding to an isotypical component $W_{{i}'}$ is the fixed subfield of the subgroup $S'\subset S$ consisting of those $g\in S$ which preserve the isotypical component: $g(W_{{i}'})=W_{{i}'}$. Let us choose one such embedding $C\rightarrow F$ (which corresponds to a choice of an isotypical component $W_{{i}'}$ of $W\otimes_k F$).\\

Let us denote by $(c_{q,g})$ the inverse matrix of the matrix $(g(e_q))$. Then $r_{\alpha,\beta}^p=\sum_{q} c_{q,Id}\cdot r_{\alpha,\beta}^{p,q}$ and for any $g\in S'=Gal(F/C)\subset S=Gal(F/k)$ if we denote by ${\phi}_g\colon A\otimes_C F\rightarrow A\otimes_C F$ the conjugation by $g\colon a\otimes f\mapsto a\otimes g(f)$, then
$$
{\phi}_g(r_{({\alpha}^i),({\beta}^i)})=g\circ r_{({\alpha}^i),({\beta}^i)}=\left( \prod_k {{\lambda}_{{\alpha}^k,{\beta}^k}}(g) \right)\cdot r_{g({\alpha}^i),g({\beta}^i)},
$$
where the action of $S'\subset S$ on upper indices $i$ (which number embeddings ${\sigma}_i\colon L\hookrightarrow F$) coincides with its action on the left cosets $S/ \tilde{H}$, where $\tilde{H}=\{ g\in S \; \mid\; g {\mid}_{{\sigma}_1(L)}=Id_{{\sigma}_1(L)} \}$ and moreover $g\in S'\subset S$ multiplies the $i$-th index ${{\alpha}}^i$ in the $r$-tuple $({\alpha}^i)=({\alpha}^1,...,{\alpha}^r)$ by $g(\theta)/\theta=\pm 1$.\\

Here ${{\lambda}_{{\alpha}^k,{\beta}^k}}(g)\in F^{*}$ are suitable constants. In order to compute them, note that isomorphisms
\begin{multline*}
{\tau}_p\colon {\wedge}_F^{p}(V\otimes_{E,{\bar{{\sigma}_i}}} F)\otimes_F (E\otimes_{E,\bar{{\sigma}_i}} F)\rightarrow {\wedge}_F^{p}(V\otimes_{E,{\bar{{\sigma}_i}}} F)\otimes_F (E\otimes_{E,\bar{{\sigma}_i}} F)\cong \\
\cong {\wedge}_F^{p}((V\otimes_{E,{{{\sigma}_i}}} F)^{*})\otimes_F (E\otimes_{E,{{\sigma}_i}} F)^{*} \rightarrow {\wedge}_F^{m-p}(V\otimes_{E,{{{\sigma}_i}}} F)\otimes_F (E\otimes_{E,{{\sigma}_i}} F)
\end{multline*}
(where the first arrow is the isomorphism determined by the matrix of ${\Phi}^{-1}$) are defined over $E$. If we assume that the isomorphism ${\wedge}_E^{p}(V)^{*} \rightarrow {\wedge}_E^{m-p}(V)\otimes_E E$ is defined via the pairing 
$$
{\wedge}_E^{p}(V) \otimes_E {\wedge}_E^{m-p}(V) \rightarrow {\wedge}_E^{m}(V)\cong E, \; x\otimes y\mapsto x\wedge y,
$$
then we find that ${{\lambda}_{{\alpha}^k,{\beta}^k}}(g)=1$, if $g(\theta)=\theta$ or ${\alpha}^k={\beta}^k$ and ${{\lambda}_{{\alpha}^k,{\beta}^k}}(g)=(-1)^{p(m-p)}\cdot (g({\sigma}_k(disc(\Phi))))^{-P({\alpha}^k,{\beta}^k)}$ otherwise.\\

Hence the matrix of $m(g)\in GL(W_{{i}'})$ is such that 
$$
m(g)\cdot E_{i,j}\cdot {m(g)}^{-1}={\phi}_g (E_{i,j})= \left( \prod_k {{\lambda}_{{\alpha}^k,{\beta}^k}}(g) \right) E_{g(i),g(j)},
$$
where $E_{i,j}$ denotes a matrix from $Mat(F)\cong End_{\Lg\otimes_k F} (W_{{i}'})$ corresponding to $r_{({\alpha}^i),({\beta}^i)}$, i.e.  conjugation by $m(g)$ acts on matrices upto a constant as the (same) permutation of columns and rows induced by $g$ on indices $({\alpha}^{i})$.\\

Then the element of $H^2(S',F^{*})$ corresponding to the central division algebra $D_0=End_{\Lg}(W_0)$ (over $C$) is the class of a $2$-cocycle $\lambda\colon S'\times S'\rightarrow F^{*}\cong F^{*}\cdot Id\subset Mat(F)$, $(g_1,g_2)\mapsto m(g_1g_2)\cdot (g_1 (m(g_2)))^{-1} \cdot m(g_1)^{-1}$ \cite{Kuznetsov}, \cite{Jacobson}.\\

\section{Example.} 

Let $k=\mathbb Q$, $r=3$ and $5\leq m\leq 6$. Let $\rho < 0$ be the negative root of the cubic polynomial $f(t)=t^3-3t+1$. Then $\frac{1}{1-\rho}$ and $1-\frac{1}{\rho}$ are the other two roots of $f(t)$ and $E=L=k(\rho)$ is a totally real cyclic cubic Galois number field \cite{Kim}.\\

Let $\Phi=-\rho\cdot X_1^2-\rho\cdot X_2^2-X_3^2-...-X_m^2$. Then by \cite{Mayanskiy} there is a $K3$ surface $X$ such that $End_{Hdg}(V)\cong E$ (where $V$ is the $\mathbb Q$-lattice of transcendental cycles on $X$), $dim_{E}V=m$ and $\Phi\colon V\otimes_E V\rightarrow E$ is the quadratic form constructed in \cite{Zarhin}.\\

Let $F=k\left( \sqrt{\rho},\sqrt{\frac{1}{1-\rho}},\sqrt{1-\frac{1}{\rho}}\right)$ be our choice of a splitting field. Note that $L\subset F$ and $\sqrt{-1}=\sqrt{\rho}\cdot\sqrt{\frac{1}{1-\rho}}\cdot\sqrt{1-\frac{1}{\rho}}\in F$. Then
$$
S=Gal(F/k)\cong ({\mathbb Z}/2{\mathbb Z})^{\oplus 3}\rtimes {\mathbb Z}/3{\mathbb Z}
$$
is a nonabelian extension of ${\mathbb Z}/3{\mathbb Z}\cong Gal(L/k)$ with generator $g$ by $({\mathbb Z}/2{\mathbb Z})^{\oplus 3}$ with generators $h_1,h_2,h_3$, where $g$ acts on the generators $(h_1,h_2,h_3)$ by the permutation $(123)$. We also denote by $g$ the element of $S$ such that $g(\sqrt{\rho})=\sqrt{\frac{1}{1-\rho}}$, $g\left(\sqrt{\frac{1}{1-\rho}}\right)=\sqrt{1-\frac{1}{\rho}}$, $g\left(\sqrt{1-\frac{1}{\rho}}\right)=\sqrt{\rho}$. We assume that each generator $h_i$, $1\leq i\leq 3$ multiplies by $-1$ the $i$-th square root among $\sqrt{\rho}, \sqrt{\frac{1}{1-\rho}}, \sqrt{1-\frac{1}{\rho}}$ and does not change the others and that $h_i {\mid}_{L}=Id$.\\

There are $3$ field embeddings $L\hookrightarrow F$: ${\sigma}_1=Id$, ${\sigma}_2=g{\mid}_L$ and ${\sigma}_3=g^2{\mid}_L$. Then $\sqrt{{{\sigma}_1}(d_1)}=\sqrt{{{\sigma}_1}(d_2)}=\sqrt{-1}\cdot \sqrt{\rho}$, $\sqrt{{{\sigma}_2}(d_1)}=\sqrt{{{\sigma}_2}(d_2)}=\sqrt{-1}\cdot \sqrt{\frac{1}{1-\rho}}$, $\sqrt{{{\sigma}_3}(d_1)}=\sqrt{{{\sigma}_3}(d_2)}=\sqrt{-1}\cdot \sqrt{1-\frac{1}{\rho}}$, $\sqrt{{{\sigma}_1}(d_3)}=\sqrt{{{\sigma}_2}(d_3)}=\sqrt{{{\sigma}_3}(d_3)}=\sqrt{-1}$, $\sqrt{-{{\sigma}_i}(d_{m-j+1})}=1$ for any $i=1,2,3$, $1\leq j\leq l=[\frac{m}{2}]$. Hence ${\otimes_{L,{\sigma}_1}} {\Gamma}_1={\otimes_{L,{\sigma}_1}} {\Gamma}_2=\sqrt{-1}\cdot \sqrt{\rho}$, ${\otimes_{L,{\sigma}_2}} {\Gamma}_1={\otimes_{L,{\sigma}_2}} {\Gamma}_2=\sqrt{-1}\cdot \sqrt{\frac{1}{1-\rho}}$, ${\otimes_{L,{\sigma}_3}} {\Gamma}_1={\otimes_{L,{\sigma}_3}} {\Gamma}_2=\sqrt{-1}\cdot \sqrt{1-\frac{1}{\rho}}$, ${\otimes_{L,{\sigma}_i}} {\Gamma}_3=\sqrt{-1}$ for all $i$ (if $m=6$).\\

(1) Let us consider first the case $m=5$. The root system is of type $B_2$: $R_0=\{ \pm {\epsilon}_p,\; \pm {\epsilon}_p \pm {\epsilon}_q \; \mid \; p,q=1,2 \}$ with basis $B_0=\{ {\epsilon}_1-{\epsilon}_2, {\epsilon}_2  \}$. Hence $B_i=\{ {\epsilon}_1 {\otimes_{L,{\sigma}_i} {\Gamma}_1}-{\epsilon}_2 {\otimes_{L,{\sigma}_i} {\Gamma}_2}, {\epsilon}_2 {\otimes_{L,{\sigma}_i} {\Gamma}_2} \}$, $1\leq i \leq 3$. The restriction of the spin representation of $\Lso(\phi)\otimes_k F$ in $C^{+}(V\otimes_k F)$ to $\Lg \otimes_k F=Res_{L/k}(\Lso(\Phi))\otimes_k F$ is isomorphic over $F$ to $2^8$ copies of the exterior tensor product ${\rho}^0 \boxtimes {\rho}^0 \boxtimes {\rho}^0$ of the irreducible spin representation of $\Lso(\Phi)\otimes_L F$. Hence over $k=\mathbb Q$ the restriction of the spin representation of $\Lso(\phi)$ in $C^{+}(V)$ to $\Lg=Res_{L/k}(\Lso(\Phi))\subset \Lso(\phi)$ is one single irreducible representation with multiplicity $\mu$ which splits over $F$ into $\frac{2^8}{\mu}$ copies of ${\rho}^0 \boxtimes {\rho}^0 \boxtimes {\rho}^0$: $C^{+}(V)\cong U^{\oplus \mu}$.\\

In order to estimate $\frac{2^8}{\mu}$ (which divides $n_0$), let us consider
$$
f_{1,...,1,1}=f_1\cdot ... \cdot f_l \cdot (1+f_0)=q\cdot \prod_{i=1}^{l} \left( e_i+\frac{\sqrt{d_i}}{\sqrt{-d_{m-i+1}}}\cdot e_{m-i+1} \right)\cdot \left( 1+\frac{1}{\sqrt{d_{l+1}}}\cdot e_{l+1} \right)
$$
(we use notation as above), where $q\in F$ is such that ${\sigma}(q)=\pm q$ for any ${\sigma}\in S=Gal(F/k)$. In our case
$$
f_{1,...,1,1}=q\cdot (e_1+\sqrt{-1}\cdot \sqrt{\rho}\cdot e_5)\cdot (e_2+\sqrt{-1}\cdot \sqrt{\rho}\cdot e_4)\cdot (1-\sqrt{-1}\cdot e_3).
$$
Hence the stabilizer of (the line in $C(V\otimes_L F)$ generated by) $f_{1,...,1,1}$ consists of the elements $g^k$, i.e. has order $3$. Since $Gal(F/k)$ has $24$ elements total, we find that $n_0=8$. Hence either $\frac{2^8}{\mu}=1$ or $\frac{2^8}{\mu}=2$ or $\frac{2^8}{\mu}=4$ or $\frac{2^8}{\mu}=8$. In the first case, ${\rho}^0 \boxtimes {\rho}^0 \boxtimes {\rho}^0$ is already defined over $\mathbb Q$ and $\mu=2^8$, while in the other cases $\mu=2^7, \mu=2^6$ and $\mu=2^5$ respectively.\\

Hence in this case $End(KS(X))_{\mathbb Q}\cong Mat_{\mu\times \mu}(D)$, where $D=End_{\Lg}(U)$ is a division algebra. Let us check that $D\cong \mathbb Q$.\\

Let us compute the cohomological invariant of $D$. In our case
$$
W\otimes_k F=V_{(1,1,1)}\oplus V_{(1,1',1')}\oplus V_{(1',1,1')}\oplus V_{(1',1',1)}\oplus V_{(2',2,2)}\oplus V_{(2,2',2)}\oplus V_{(2,2,2')}\oplus V_{(2',2',2')},
$$
where $V_{(p_1,p_2,p_3)}=S_{p_1}^1\otimes_F S_{p_2}^2\otimes_F S_{p_3}^3$ in the notation of Section 5.2 and the values $1,1',2,2'$ of $p_i$ correspond to the indices $({\alpha}_1,{\alpha}_2,\gamma)$ of ideals $I_{{\alpha}_1,{\alpha}_2,\gamma}$ as follows: $1=(+++)$, $1'=(--+)$, $2=(---)$, $2'=(++-)$.\\

Let us denote $\bar{1}=2$, $\bar{1'}=2'$, $\bar{2}=1$, $\bar{2'}=1'$ and $\tilde{1}=2'$, $\tilde{1'}=2$, $\tilde{2}=1'$, $\tilde{2'}=1$. Then $g(V_{(p_1,p_2,p_3)})=V_{(p_3,p_1,p_2)}$ and $h_i(V_{(p_1,p_2,p_3)})= V_{(q_1,q_2,q_3)}$, where $q_i=\tilde{p_i}$ and $q_j=\bar{p_j}$ for $j\neq i$.\\

Let us denote $a=(1,1',1')$, $b=(1',1,1')$, $c=(1',1',1)$, $d=(1,1,1)$, $p=(2',2,2)$, $q=(2,2',2)$, $r=(2,2,2')$, $s=(2',2',2')$. Then using formulas from Section 6 we can choose coefficients ${\lambda}_{\alpha,\beta}=\prod_{i=1}^r {\lambda}_{{\alpha}^i,{\beta}^i}\in F^{*}$ as follows:
\begin{itemize}
\item ${\lambda}_{\alpha,\beta}=1$ for $(\alpha,\beta)\in \{ (d,-), (s,-), (a,a), (b,b), (c,c), (p,p), (q,q), (r,r) \}$, 
\item ${\lambda}_{\alpha,\beta}=1$ for $(\alpha,\beta)\in \{ (b,q), (a,p), (c,r), (q,b), (p,a), (r,c) \}$, 
\item ${\lambda}_{\alpha,\beta}=c_1$ for $(\alpha,\beta)\in \{ (b,a), (b,p), (c,a), (c,p), (q,a), (q,p), (r,a), (r,p) \}$, 
\item ${\lambda}_{\alpha,\beta}=c_2$ for $(\alpha,\beta)\in \{ (a,b), (a,q), (c,b), (c,q), (p,b), (p,q), (r,b), (r,q) \}$,
\item ${\lambda}_{\alpha,\beta}=c_3$ for $(\alpha,\beta)\in \{ (b,c), (b,r), (a,c), (a,r), (p,c), (p,r), (q,c), (q,r) \}$,
\item ${\lambda}_{\alpha,\beta}=c_1c_2$ for $(\alpha,\beta)\in \{ (c,d), (c,s), (r,d), (r,s) \}$,
\item ${\lambda}_{\alpha,\beta}=c_1c_3$ for $(\alpha,\beta)\in \{ (b,d), (b,s), (q,d), (q,s) \}$,
\item ${\lambda}_{\alpha,\beta}=c_2c_3$ for $(\alpha,\beta)\in \{ (a,d), (a,s), (p,d), (p,s) \}$.
\end{itemize}

Here we denoted $c_i={\sigma}_i\left(\frac{-1}{\Phi(f_1,f_{-1})\cdot \Phi(f_2,f_{-2})}\right)={\sigma}_i\left(\frac{-1}{4{\rho}^2}\right)$.\\

Then in the formulas in Section 6 we can take:

\begin{itemize}
\item $m(g)= \begin{pmatrix} G & 0\\ 0 & G \end{pmatrix}$ is an $8\times 8$ matrix whose rows and columns are numbered according to the following sequence of indices of $V_{(p_1,p_2,p_3)}: \; (dabcspqr)$,  
\item $m(h_1)= \begin{pmatrix} 0 & X_1^{-1}\\ \frac{1}{c_2c_3}\cdot X_1 & 0 \end{pmatrix}$ is an $8\times 8$ matrix whose rows and columns are numbered according to the following sequence of indices of $V_{(p_1,p_2,p_3)}: \; (dabcpsrq)$,
\item $m(h_2)= \begin{pmatrix} 0 & X_2^{-1}\\ \frac{1}{c_1c_3}\cdot X_2 & 0 \end{pmatrix}$ is an $8\times 8$ matrix whose rows and columns are numbered according to the following sequence of indices of $V_{(p_1,p_2,p_3)}: \; (dabcqrsp)$,
\item $m(h_3)= \begin{pmatrix} 0 & X_3^{-1}\\ \frac{1}{c_1c_2}\cdot X_3 & 0 \end{pmatrix}$ is an $8\times 8$ matrix whose rows and columns are numbered according to the following sequence of indices of $V_{(p_1,p_2,p_3)}: \; (dabcrqps)$,
\item $m(g^k\cdot h_1^{a_1} h_2^{a_2} h_3^{a_3})=m(g)^k\cdot g^k \left( m(h_1)^{a_1}\cdot m(h_2)^{a_2}\cdot m(h_3)^{a_3} \right)$, where $0\leq a_i\leq 1$, $k\geq 0$. 
\end{itemize}

Here we denoted $G=\begin{pmatrix} 1 & 0 & 0 & 0\\ 0 & 0 & 0 & 1\\ 0 & 1 & 0 & 0\\ 0 & 0 & 1 & 0 \end{pmatrix}$, $X_1=\begin{pmatrix} 1 & 0 & 0 & 0\\ 0 & c_2c_3 & 0 & 0\\ 0 & 0 & c_3 & 0\\ 0 & 0 & 0 & c_2 \end{pmatrix}$, $X_2=\begin{pmatrix} 1 & 0 & 0 & 0\\ 0 & c_3 & 0 & 0\\ 0 & 0 & c_1c_3 & 0\\ 0 & 0 & 0 & c_1 \end{pmatrix}$ and $X_3=\begin{pmatrix} 1 & 0 & 0 & 0\\ 0 & c_2 & 0 & 0\\ 0 & 0 & c_1 & 0\\ 0 & 0 & 0 & c_1c_2 \end{pmatrix}$.\\

Note that $m(h_i)\cdot m(h_j)=m(h_j)\cdot m(h_i)$, $m(h_i)^2=\frac{c_i}{c_1c_2c_3}$, $m(g)^3=1$ and $m(gh_ig^{-1})=m(g)\cdot g(m(h_i))\cdot m(g)^{-1}$.\\

This implies that the class of $D$ in $H^2(S,F^{*})$ is represented by the $2$-cocycle $\lambda \colon S\times S\rightarrow F^{*}$ such that ${\lambda}(h_1^{a_1}h_2^{a_2}h_3^{a_3},h_1^{b_1}h_2^{b_2}h_3^{b_3})=(c_2c_3)^{x_1}\cdot (c_1c_3)^{x_2}\cdot (c_1c_2)^{x_3}$ and ${\lambda}(g^kh,g^lh')=g^{k+l}({\lambda}(g^{-l}hg^l,h'))$, where $0\leq a_i\leq 1$, $0\leq b_i\leq 1$, $x_i=1$ if $a_i=b_i=1$ and $0$ otherwise, and $h, h'$ are elements of the subgroup $({\mathbb Z}/2{\mathbb Z})^{\oplus 3}\subset S$ generated by $h_1, h_2, h_3$.\\

Since $c_ic_j=\left( \frac{1}{4\cdot {\sigma}_i(\rho){\sigma}_j(\rho)} \right)^2$ is a square in $L^{*}$, we conclude that $\lambda$ is a coboundary. Namely, the required morphism $c\colon S\rightarrow F^{*}$ (whose coboundary is $\lambda$) can be defined as follows:
$$
c(g^k\cdot h_1^{a_1} h_2^{a_2} h_3^{a_3})=g^k\left( (\sqrt{c_2c_3})^{a_1}\cdot (\sqrt{c_1c_3})^{a_2}\cdot  (\sqrt{c_1c_2})^{a_3} \right),
$$
where $0\leq a_i\leq 1$, $k\geq 0$. Note that $c(gh_ig^{-1})=g(c(h_i))$. So, the class of $D$ in $H^2(S,F^{*})$ vanishes. Hence $D\cong \mathbb Q$.\\ 

So, in this example $End(KS(X))_{\mathbb Q}\cong Mat_{256\times 256}(\mathbb Q)$.\\

(2) Now let us consider the case $m=6$. The root system is of type $D_3$: $R_0=\{ \pm {\epsilon}_p \pm {\epsilon}_q \; \mid \; p,q=1,2,3 \}$ with basis $B_0=\{ {\epsilon}_1-{\epsilon}_2, {\epsilon}_2-{\epsilon}_3, {\epsilon}_2+{\epsilon}_3  \}$. Hence $B_i=\{ {\epsilon}_1 {\otimes_{L,{\sigma}_i} {\Gamma}_1}-{\epsilon}_2 {\otimes_{L,{\sigma}_i} {\Gamma}_2}, {\epsilon}_2 {\otimes_{L,{\sigma}_i} {\Gamma}_2}-{\epsilon}_3 {\otimes_{L,{\sigma}_i} {\Gamma}_3}, {\epsilon}_2 {\otimes_{L,{\sigma}_i} {\Gamma}_2}+{\epsilon}_3 {\otimes_{L,{\sigma}_i} {\Gamma}_3} \}$, $1\leq i \leq 3$, and the Weyl group is generated by sign inversions in front of two of ${\epsilon}_1, {\epsilon}_2, {\epsilon}_3$ and by all possible permutations of ${\epsilon}_1, {\epsilon}_2, {\epsilon}_3$.\\

The restriction of the spin representation of $\Lso(\phi)\otimes_k F$ in $C^{+}(V\otimes_k F)$ to $\Lg \otimes_k F=Res_{L/k}(\Lso(\Phi))\otimes_k F$ is isomorphic over $F$ to the sum of the exterior tensor products of semi-spin representations (in all possible combinations) each with multiplicity $2^8$: $C^{+}(V\otimes_k F)\cong \bigoplus_{{\alpha}_1,{\alpha}_2,{\alpha}_3 \in \{ \pm \} } 2^8\cdot ({\rho}^{{\alpha}_1} \boxtimes {\rho}^{{\alpha}_2} \boxtimes {\rho}^{{\alpha}_3})$. Hence the set $\Omega$ of highest weights consists of the elements ${\omega}_{{\alpha}_1,{\alpha}_2,{\alpha}_3}=\frac{1}{2}\cdot \sum_{i=1}^{3}( {\epsilon}_1 {\otimes_{L,{\sigma}_i} {\Gamma}_1}+ {\epsilon}_2 {\otimes_{L,{\sigma}_i} {\Gamma}_2} + {\alpha}_i\cdot {\epsilon}_3 {\otimes_{L,{\sigma}_i} {\Gamma}_3})$ for various ${\alpha}_i\in \{ \pm 1  \}$.\\

Note that $g({\omega}_{{\alpha}_1,{\alpha}_2,{\alpha}_3})={\omega}_{{\alpha}_3,{\alpha}_1,{\alpha}_2}$ and $h_i({\omega}_{{\alpha}_1,{\alpha}_2,{\alpha}_3})={\omega}_{-{\alpha}_1,-{\alpha}_2,-{\alpha}_3}$. So, $\Omega={\Omega}_1\cup {\Omega}_2$ has two $S$-orbits: ${\Omega}_1=\{ {\omega}_{+,+,+},{\omega}_{-,-,-}  \}$ and ${\Omega}_2=\{ {\omega}_{+,+,-},{\omega}_{+,-,+},{\omega}_{-,+,+},{\omega}_{-,-,+}, {\omega}_{-,+,-} , {\omega}_{+,-,-} \}$.\\

Hence over $k=\mathbb Q$ we have: $C^{+}(V)\cong U^{\oplus \mu} \oplus V^{\oplus \nu}$ as $\Lg$-modules, where $U$ and $V$ are not isomorphic as representations of $\Lg=Res_{L/k}(\Lso(\Phi))$. $U\otimes_k F$ splits into $\frac{2^8}{\mu}$ copies of ${\rho}^{+} \boxtimes {\rho}^{+} \boxtimes {\rho}^{+}$ and $\frac{2^8}{\mu}$ copies of ${\rho}^{-} \boxtimes {\rho}^{-} \boxtimes {\rho}^{-}$, while $V\otimes_k F$ splits into $\frac{2^8}{\nu}$ copies of ${\rho}^{{\alpha}_1} \boxtimes {\rho}^{{\alpha}_2} \boxtimes {\rho}^{{\alpha}_3}$ with other ${{\alpha}_i}$'s.\\

In order to estimate multiplicities $\mu$ and $\nu$, let us consider
$$
f_{1,...,1}=f_1\cdot ... \cdot f_l=q\cdot \prod_{i=1}^{l} \left( e_i+\frac{\sqrt{d_i}}{\sqrt{-d_{m-i+1}}}\cdot e_{m-i+1} \right)
$$
(we use notation introduced above), where $q\in F$ is such that ${\sigma}(q)=\pm q$ for any ${\sigma}\in S=Gal(F/k)$. In our case
$$
f_{1,...,1}=q\cdot (e_1+\sqrt{-1}\cdot \sqrt{\rho}\cdot e_6)\cdot (e_2+\sqrt{-1}\cdot \sqrt{\rho}\cdot e_5)\cdot (e_3+\sqrt{-1} \cdot e_4).
$$
Hence the stabilizer of (the line in $C(V\otimes_L F)$ generated by) $f_{1,...,1}$ consists of the elements $g^k$. Since the stabilizer of ${\omega}_{+,+,+}\in\Omega$ as a subgroup of $S$ is generated by elements $g,h_1h_2,h_1h_3,h_2h_3$, we conclude that $n_{+,+,+}=4$. Since the stabilizer of ${\omega}_{+,+,-}\in\Omega$ has $4$ elements: $Id$ and $h_1h_2, h_1h_3, h_2h_3$, we conclude that $n_{+,+,-}=4$ as well. The same computation as in the case $m=5$ above shows that $\frac{2^8}{\mu}=\frac{2^8}{\nu}=1$, i.e. $\mu=\nu=256$, and the division algebras $D_1=End_{\Lg}(U)$ and $D_2=End_{\Lg}(V)$ are fields, i.e. coincide with their centers.\\

According to Section 6.2, the center $C_1$ of $D_1$ is the subfield of $F$ fixed by the stabilizer of ${\omega}_{+,+,+}\in\Omega$, i.e. $D_1=C_1\cong k(\sqrt{-1})$. Similarly, the center $C_2$ of $D_2$ is the subfield of $F$ fixed by the stabilizer of ${\omega}_{+,+,-}\in\Omega$, i.e. $D_2=C_2\cong k(\sqrt{-1},\rho)$.\\

So, in this example $End(KS(X))_{\mathbb Q}\cong Mat_{256\times 256}(\mathbb Q (\sqrt{-1})) \times Mat_{256\times 256}(\mathbb Q (\sqrt{-1},\rho) )$.\\

\medskip

(3) Let us modify the first example above. Consider the same number $\rho$ and the same totally real cubic field $E=L=k(\rho)$, but a different quadratic form
$$
\Phi=-(a+\rho)\cdot X_1^2-(a+\rho)\cdot X_2^2-X_3^2-X_4^2-X_5^2,
$$
where $a$ is a fixed rational number between $0$ and $-\rho$: $0<a<-\rho$. As above, these quadratic form and totally real field correspond to a $K3$ surface $X$ (\cite{Mayanskiy}). Assume that $1+3a-a^3>0$ is not a square of a rational number.\\

Let $F=k\left(\sqrt{a+\rho},\sqrt{a+\frac{1}{1-\rho}},\sqrt{a+1-\frac{1}{\rho}},\sqrt{-1}\right)$ be our choice of a splitting field. Note that $L\subset F$ and $\sqrt{a+\rho}\cdot \sqrt{a+\frac{1}{1-\rho}}\cdot \sqrt{a+1-\frac{1}{\rho}}=\sqrt{-1-3a+a^3}$. Then 
$$
S=Gal(F/k)\cong {\mathbb Z}/2{\mathbb Z}\oplus G,
$$
where $G$ is the group isomorphic to the Galois group of the splitting field from the first example above, i.e. $G$ is a noncommutative group extension of ${\mathbb Z}/3{\mathbb Z}\cong Gal(L/k)$ by $({\mathbb Z}/2{\mathbb Z})^{\oplus 3}$. Let $g$ be a generator of ${\mathbb Z}/3{\mathbb Z}$ such that $g(\sqrt{a+\rho})=\sqrt{a+\frac{1}{1-\rho}}$, $g\left(\sqrt{a+\frac{1}{1-\rho}}\right)=\sqrt{a+1-\frac{1}{\rho}}$, $g\left(\sqrt{a+1-\frac{1}{\rho}}\right)=\sqrt{a+\rho}$, $g(\sqrt{-1})=\sqrt{-1}$. Let $h_1,h_2,h_3$ be the generators of $({\mathbb Z}/2{\mathbb Z})^{\oplus 3}$ and $h_0$ be the generator of the first factor ${\mathbb Z}/2{\mathbb Z}$ in $S$ above such that each $h_i$, $0\leq i\leq 3$ multiplies by $-1$ the $i$-th square root among $\sqrt{-1}, \sqrt{a+\rho}, \sqrt{a+\frac{1}{1-\rho}}, \sqrt{a+1-\frac{1}{\rho}}$ and does not change the others. We also assume that $h_i {\mid}_{L}=Id$, $0\leq i\leq 3$.\\

There are $3$ field embeddings $L\hookrightarrow F$: ${\sigma}_1=Id$, ${\sigma}_2=g{\mid}_L$ and ${\sigma}_3=g^2{\mid}_L$. Then $\sqrt{{{\sigma}_1}(d_1)}=\sqrt{{{\sigma}_1}(d_2)}=\sqrt{-1}\cdot \sqrt{a+\rho}$, $\sqrt{{{\sigma}_2}(d_1)}=\sqrt{{{\sigma}_2}(d_2)}=\sqrt{-1}\cdot \sqrt{a+\frac{1}{1-\rho}}$, $\sqrt{{{\sigma}_3}(d_1)}=\sqrt{{{\sigma}_3}(d_2)}=\sqrt{-1}\cdot \sqrt{a+1-\frac{1}{\rho}}$, $\sqrt{{{\sigma}_1}(d_3)}=\sqrt{{{\sigma}_2}(d_3)}=\sqrt{{{\sigma}_3}(d_3)}=\sqrt{-1}$, $\sqrt{-{{\sigma}_i}(d_4)}=\sqrt{-{{\sigma}_i}(d_5)}=1$ for any $i=1,2,3$. Hence ${\otimes_{L,{\sigma}_1}} {\Gamma}_1={\otimes_{L,{\sigma}_1}} {\Gamma}_2=\sqrt{-1}\cdot \sqrt{a+\rho}$, ${\otimes_{L,{\sigma}_2}} {\Gamma}_1={\otimes_{L,{\sigma}_2}} {\Gamma}_2=\sqrt{-1}\cdot \sqrt{a+\frac{1}{1-\rho}}$, ${\otimes_{L,{\sigma}_3}} {\Gamma}_1={\otimes_{L,{\sigma}_3}} {\Gamma}_2=\sqrt{-1}\cdot \sqrt{a+1-\frac{1}{\rho}}$.\\

As in the first example above, the root system is of type $B_2$: $R_0=\{ \pm {\epsilon}_p,\; \pm {\epsilon}_p \pm {\epsilon}_q \; \mid \; p,q=1,2 \}$ with basis $B_0=\{ {\epsilon}_1-{\epsilon}_2, {\epsilon}_2  \}$. Hence $B_i=\{ {\epsilon}_1 {\otimes_{L,{\sigma}_i} {\Gamma}_1}-{\epsilon}_2 {\otimes_{L,{\sigma}_i} {\Gamma}_2}, {\epsilon}_2 {\otimes_{L,{\sigma}_i} {\Gamma}_2} \}$, $1\leq i \leq 3$. The restriction of the spin representation of $\Lso(\phi)\otimes_k F$ in $C^{+}(V\otimes_k F)$ to $\Lg \otimes_k F=Res_{L/k}(\Lso(\Phi))\otimes_k F$ is isomorphic over $F$ to $2^8$ copies of the exterior tensor product ${\rho}^0 \boxtimes {\rho}^0 \boxtimes {\rho}^0$ of the irreducible spin representation of $\Lso(\Phi)\otimes_L F$. Hence over $k=\mathbb Q$ the restriction of the spin representation of $\Lso(\phi)$ in $C^{+}(V)$ to $\Lg=Res_{L/k}(\Lso(\Phi))\subset \Lso(\phi)$ is one single irreducible representation with multiplicity $\mu$ which splits over $F$ into $\frac{2^8}{\mu}$ copies of ${\rho}^0 \boxtimes {\rho}^0 \boxtimes {\rho}^0$: $C^{+}(V)\cong U^{\oplus \mu}$.\\

In order to estimate $\frac{2^8}{\mu}$ (which divides $n_0$), let us consider
$$
f_{1,...,1,1}=f_1\cdot ... \cdot f_l \cdot (1+f_0)=q\cdot \prod_{i=1}^{l} \left( e_i+\frac{\sqrt{d_i}}{\sqrt{-d_{m-i+1}}}\cdot e_{m-i+1} \right)\cdot \left( 1+\frac{1}{\sqrt{d_{l+1}}}\cdot e_{l+1} \right)
$$
(we use notation as above), where $q\in F$ is such that ${\sigma}(q)=\pm q$ for any ${\sigma}\in S=Gal(F/k)$. In our case
$$
f_{1,...,1,1}=q\cdot (e_1+\sqrt{-1}\cdot \sqrt{a+\rho}\cdot e_5)\cdot (e_2+\sqrt{-1}\cdot \sqrt{a+\rho}\cdot e_4)\cdot (1-\sqrt{-1}\cdot e_3).
$$
Hence the stabilizer of (the line in $C(V\otimes_L F)$ generated by) $f_{1,...,1,1}$ consists of the elements $g^k$, i.e. has order $3$. Since $Gal(F/k)$ has $48$ elements total, we find that $n_0=16$. Hence either $\frac{2^8}{\mu}=1$ or $\frac{2^8}{\mu}=2$ or $\frac{2^8}{\mu}=4$ or $\frac{2^8}{\mu}=8$ or $\frac{2^8}{\mu}=16$. In the first case, ${\rho}^0 \boxtimes {\rho}^0 \boxtimes {\rho}^0$ is already defined over $\mathbb Q$ and $\mu=2^8$, while in the other cases $\mu=2^7, \mu=2^6$, $\mu=2^5$ and $\mu=2^4$ respectively.\\

Hence in this case $End(KS(X))_{\mathbb Q}\cong Mat_{\mu\times \mu}(D)$, where $D=End_{\Lg}(U)$ is a division algebra. Let us compute the cohomological invariant of $D$. In our case 
\begin{multline*}
W\otimes_k F=V_{(1,1,1)}\oplus V_{(1',1,1)}\oplus V_{(1,1',1)}\oplus V_{(1,1,1')}\oplus V_{(1',1',1')}\oplus V_{(1,1',1')}\oplus V_{(1',1,1')}\oplus V_{(1',1',1)}\oplus \\
\oplus V_{(2,2,2)}\oplus V_{(2',2,2)}\oplus V_{(2,2',2)}\oplus V_{(2,2,2')}\oplus V_{(2',2',2')}\oplus V_{(2,2',2')}\oplus V_{(2',2,2')}\oplus V_{(2',2',2)},
\end{multline*}
where $V_{(p_1,p_2,p_3)}=S_{p_1}^1\otimes_F S_{p_2}^2\otimes_F S_{p_3}^3$ in the notation of Section 5.2 and the values $1,1',2,2'$ of $p_i$ correspond to the indices $({\alpha}_1,{\alpha}_2,\gamma)$ of ideals $I_{{\alpha}_1,{\alpha}_2,\gamma}$ as follows: $1=(+++)$, $1'=(--+)$, $2=(++-)$, $2'=(---)$.\\

Let us denote $\bar{1}=1'$, $\bar{1'}=1$, $\bar{2}=2'$, $\bar{2'}=2$ and $\tilde{1}=2'$, $\tilde{1'}=2$, $\tilde{2}=1'$, $\tilde{2'}=1$. Then $g(V_{(p_1,p_2,p_3)})=V_{(p_3,p_1,p_2)}$, $h_0(V_{(p_1,p_2,p_3)})= V_{(\tilde{p_1},\tilde{p_2},\tilde{p_3})}$ and $h_i(V_{(p_1,p_2,p_3)})= V_{(q_1,q_2,q_3)}$, $1\leq i\leq 3$, where $q_i=\bar{p_i}$ and $q_j={p_j}$ for $j\neq i$.\\

Let us denote $a=(1',1,1)$, $b=(1,1',1)$, $c=(1,1,1')$, $d=(1,1,1)$, $a'=(1,1',1')$, $b'=(1',1,1')$, $c'=(1',1',1)$, $d'=(1',1',1')$, $p=(2',2,2)$, $q=(2,2',2)$, $r=(2,2,2')$, $s=(2,2,2)$, $p'=(2,2',2')$, $q'=(2',2,2')$, $r'=(2',2',2)$, $s'=(2',2',2')$. Consider the set of indices $T=\{ d,a,b,c,d',a',b',c',s,p,q,r,s',p',q',r'  \}$ and the morphism $t\colon T\rightarrow T, s\mapsto d, p\mapsto a, q\mapsto b, r\mapsto c, s'\mapsto d', p'\mapsto a', q'\mapsto b', r'\mapsto c'$ and $x\mapsto x$ for all other $x\in T$.\\

Then using formulas from Section 6 we can choose coefficients ${\lambda}_{\alpha,\beta}=\prod_{i=1}^r {\lambda}_{{\alpha}^i,{\beta}^i}\in F^{*}$ as follows:
\begin{itemize}
\item ${\lambda}_{d,x}={\lambda}_{s,x}={\lambda}_{x,x}=1$, ${\lambda}_{d',d}=c_1c_2c_3$ and ${\lambda}_{x,y}={\lambda}_{t(x),t(y)}$ for any $x,y\in T$,
\item ${\lambda}_{\alpha,\beta}=1$ for $(\alpha,\beta)\in \{ (a,d'), (a,b'), (a,c'), (b,d'), (b,a'), (b,c') \}$, 
\item ${\lambda}_{\alpha,\beta}=1$ for $(\alpha,\beta)\in \{ (c,d'), (c,b'), (c,a'), (a',d'), (b',d'), (c',d')  \}$, 
\item ${\lambda}_{\alpha,\beta}=c_1$ for $(\alpha,\beta)\in \{ (a,a'), (a,d), (a,b), (a,c), (d',a'), (b',c), (b',a'), (c',b), (c',a') \}$, 
\item ${\lambda}_{\alpha,\beta}=c_2$ for $(\alpha,\beta)\in \{ (b,b'), (b,d), (b,a), (b,c), (d',b'), (a',c), (a',b'), (c',a), (c',b') \}$,
\item ${\lambda}_{\alpha,\beta}=c_3$ for $(\alpha,\beta)\in \{ (c,c'), (c,d), (c,b), (c,a), (d',c'), (a',b), (a',c'), (b',a), (b',c') \}$,
\item ${\lambda}_{\alpha,\beta}=c_1c_2$ for $(\alpha,\beta)\in \{ (d',c), (c',c), (c',d) \}$,
\item ${\lambda}_{\alpha,\beta}=c_1c_3$ for $(\alpha,\beta)\in \{ (d',b), (b',b), (b',d) \}$,
\item ${\lambda}_{\alpha,\beta}=c_2c_3$ for $(\alpha,\beta)\in \{ (d',a), (a',a), (a',d) \}$.
\end{itemize}

Here we denoted $c_i={\sigma}_i\left(\frac{-1}{\Phi(f_1,f_{-1})\cdot \Phi(f_2,f_{-2})}\right)={\sigma}_i\left(\frac{-1}{4(a+{\rho})^2}\right)$.\\

Then in the formulas in Section 6 we can take:
\begin{itemize}
\item $m(g)= \begin{pmatrix} G & 0 & 0 & 0\\ 0 & G & 0 & 0\\ 0 & 0 & G & 0\\ 0 & 0 & 0 & G \end{pmatrix}$ is a $16\times 16$ matrix whose rows and columns are numbered according to the following sequence of indices of $V_{(p_1,p_2,p_3)}: \; (dabcd'a'b'c'spqrs'p'q'r')$,  
\item $m(h_1)= \begin{pmatrix} 0 & 1 & 0 & 0\\ \frac{1}{c_1}\cdot 1 & 0 & 0 & 0\\ 0 & 0 & 0 & 1\\ 0 & 0 & \frac{1}{c_1}\cdot 1 & 0 \end{pmatrix}$ is a $16\times 16$ matrix whose rows and columns are numbered according to the following sequence of indices of $V_{(p_1,p_2,p_3)}: \; (da'bcad'c'b'sp'qrps'r'q')$,
\item $m(h_2)= \begin{pmatrix} 0 & 1 & 0 & 0\\ \frac{1}{c_2}\cdot 1 & 0 & 0 & 0\\ 0 & 0 & 0 & 1\\ 0 & 0 & \frac{1}{c_2}\cdot 1 & 0 \end{pmatrix}$ is a $16\times 16$ matrix whose rows and columns are numbered according to the following sequence of indices of $V_{(p_1,p_2,p_3)}: \; (dab'cbc'd'a'spq'rqr's'p')$,
\item $m(h_3)= \begin{pmatrix} 0 & 1 & 0 & 0\\ \frac{1}{c_3}\cdot 1 & 0 & 0 & 0\\ 0 & 0 & 0 & 1\\ 0 & 0 & \frac{1}{c_3}\cdot 1 & 0 \end{pmatrix}$ is a $16\times 16$ matrix whose rows and columns are numbered according to the following sequence of indices of $V_{(p_1,p_2,p_3)}: \; (dabc'cb'a'd'spqr'rq'p's')$,
\item $m(h_0)= \begin{pmatrix} 0 & X_0^{-1} & 0 & 0\\ \frac{1}{c_1c_2c_3}\cdot X_0 & 0 & 0 & 0\\ 0 & 0 & 0 & \frac{1}{c_1c_2c_3}\cdot X_0\\ 0 & 0 & X_0^{-1} & 0 \end{pmatrix}$ is a $16\times 16$ matrix whose rows and columns are numbered according to the following sequence of indices of $V_{(p_1,p_2,p_3)}: \; (dabcs'p'q'r'd'a'b'c'spqr)$,
\item $m(g^k\cdot h_0^{a_0} h_1^{a_1} h_2^{a_2} h_3^{a_3})=m(g)^k\cdot g^k \left( m(h_0)^{a_0}\cdot m(h_1)^{a_1}\cdot m(h_2)^{a_2}\cdot m(h_3)^{a_3} \right)$, where $0\leq a_i\leq 1$, $k\geq 0$. 
\end{itemize}

Here we denoted $G=\begin{pmatrix} 1 & 0 & 0 & 0\\ 0 & 0 & 0 & 1\\ 0 & 1 & 0 & 0\\ 0 & 0 & 1 & 0 \end{pmatrix}$, $1=\begin{pmatrix} 1 & 0 & 0 & 0\\ 0 & 1 & 0 & 0\\ 0 & 0 & 1 & 0\\ 0 & 0 & 0 & 1 \end{pmatrix}$ (in the definitions of $m(h_i)$) and $X_0=\begin{pmatrix} 1 & 0 & 0 & 0\\ 0 & c_1 & 0 & 0\\ 0 & 0 & c_2 & 0\\ 0 & 0 & 0 & c_3 \end{pmatrix}$.\\

Note that $m(h_i)\cdot m(h_j)=m(h_j)\cdot m(h_i)$, $m(h_i)^2=\frac{1}{c_i}$, $1\leq i\leq 3$, $m(h_0)^2=\frac{1}{c_1c_2c_3}$, $m(g)^3=1$ and $m(gh_ig^{-1})=m(g)\cdot g(m(h_i))\cdot m(g)^{-1}$.\\

This implies that the class of $D$ in $H^2(S,F^{*})$ is represented by the $2$-cocycle $\lambda \colon S\times S\rightarrow F^{*}$ such that
$$
{\lambda}(h_0^{a_0}h_1^{a_1}h_2^{a_2}h_3^{a_3},h_0^{b_0}h_1^{b_1}h_2^{b_2}h_3^{b_3})=(c_1c_2c_3)^{x_0}\cdot (c_2c_3)^{x_1}\cdot (c_1c_3)^{x_2}\cdot (c_1c_2)^{x_3}
$$
and ${\lambda}(g^kh,g^lh')=g^{k+l}({\lambda}(g^{-l}hg^l,h'))$, where $0\leq a_i\leq 1$, $0\leq b_i\leq 1$, $x_i=1$ if $a_i=b_i=1$ and $0$ otherwise, and $h, h'$ are elements of the subgroup ${\mathbb Z}/2{\mathbb Z}\oplus({\mathbb Z}/2{\mathbb Z})^{\oplus 3}\subset S$ generated by $h_0,h_1, h_2, h_3$.\\

Let us multiply $\lambda$ by the inverse of the coboundary of the 1-cochain given by the morphism $c\colon S\rightarrow F^{*}$ such that
$$
c(g^k\cdot h_0^{a_0} h_1^{a_1} h_2^{a_2} h_3^{a_3})=g^k\left( (\sqrt{c_1c_2c_3})^{a_0}\cdot (\sqrt{c_1})^{a_1}\cdot (\sqrt{c_2})^{a_2}\cdot  (\sqrt{c_3})^{a_3} \right),
$$
where $0\leq a_i\leq 1$, $k\geq 0$. Note that $c(gh_ig^{-1})=g(c(h_i))$.\\

This changes $\lambda$ to a 2-cocycle ${\lambda}' \colon S\times S\rightarrow F^{*}$ such that
$$
{\lambda}'(g^k\cdot h_0^{a_0}h_1^{a_1}h_2^{a_2}h_3^{a_3},g^l\cdot h_0^{b_0}h_1^{b_1}h_2^{b_2}h_3^{b_3})=(-1)^{a_0\cdot (b_0+b_1+b_2+b_3)},
$$
where $0\leq a_i\leq 1$, $0\leq b_i\leq 1$.\\

Let $H\subset S$ be the subgroup generated by $g, h_1h_2, h_1h_3, h_2h_3$ and
$$
F^H=k\left(\sqrt{-1},\sqrt{(a+\rho)(a+\frac{1}{1-\rho})(a+1-\frac{1}{\rho})} \right)=k(\sqrt{-1},\sqrt{-1-3a+a^3})
$$
be the corresponding fixed subfield of $F$. Denote the generators of $Gal(F^H/k)\cong {\mathbb Z}/2{\mathbb Z}\oplus {\mathbb Z}/2{\mathbb Z}$ by $h_0$ and $h=h_1h_2h_3$.\\

We see that the class of $D$ in $H^2(S,F^{*})$ is the image under the inflation homomorphism $H^2(Gal(F^H/k),(F^H)^{*})\rightarrow H^2(S,F^{*})$ of a class represented by the $2$-cocycle ${\lambda}''\colon Gal(F^H/k)\times Gal(F^H/k) \rightarrow k(\sqrt{-1},\sqrt{-1-3a+a^3})^{*}$ such that ${\lambda}''(h_0,h_0)={\lambda}''(h_0,h)=-{\lambda}''(h,h_0)=-{\lambda}''(h,h)=-1$. Multiplying it by the coboundary of the 1-cochain given by the morphism $c\colon Gal(F^H/k)\rightarrow (F^H)^{*}$ such that $c(h)=c(h_0)=\sqrt{-1}$, $c(hh_0)=1$, we obtain a 2-cocycle (also denoted by ${\lambda}''$) with the property ${\lambda}''(h_0h,-)={\lambda}''(-,h_0h)=1$ and ${\lambda}''(h_0,h_0)=1$. Note that $(F^H)^{<h_0h>}=k(\sqrt{1+3a-a^3})$ is a totally real quadratic field with Galois group ${\mathbb Z}/2{\mathbb Z}$ with generator $1$.\\

This means that the cohomological class of $D$ can be obtained via the inflation homomorphism from the class in $H^2(Gal(k(\sqrt{1+3a-a^3})/k), k(\sqrt{1+3a-a^3})^{*})$ of the $2$-cocycle ${\lambda}_0\colon {\mathbb Z}/2{\mathbb Z}\times {\mathbb Z}/2{\mathbb Z} \rightarrow k(\sqrt{1+3a-a^3})^{*}$ such that ${\lambda}_0(1,1)=-1$.\\

Hence $D$ is a quaternion algebra over $\mathbb Q=k$ of degree $deg(D)=2$ split over ${\mathbb Q}(\sqrt{1+3a-a^3})$ with $4$ generators over $\mathbb Q$: $1,i,j,k$ such that $i^2=j^2=1+3a-a^3$, $k=ij=-ji$. In other words, $D=(1+3a-a^3,1+3a-a^3)_{\mathbb Q}$.\\

So, in this example $End(KS(X))_{\mathbb Q}\cong Mat_{128\times 128}((1+3a-a^3,1+3a-a^3)_{\mathbb Q})$.\\

\medskip

(4) If in the previous example we take
$$
\Phi=-(b\cdot\rho)\cdot X_1^2-(b\cdot\rho)\cdot X_2^2-X_3^2-X_4^2-X_5^2,
$$ 
where $b> 0$ is a rational number which is not a square of another rational number, then the same computation as above gives:
$$
End(KS(X))_{\mathbb Q}\cong Mat_{128\times 128}((b,b)_{\mathbb Q})
$$
for the corresponding $K3$ surface $X$.\\

\section{Acknowledgement.} 

We thank Yuri Zarhin for suggesting this problem and for pointing out an error in the Example section of the previous version. Many of our constructions were influenced by papers \cite{vanGeemen}, \cite{vanGeemen1} and \cite{vanGeemen2}, where Bert van Geemen studies endomorphism algebras of Kuga-Satake varieties in more special cases.\\

\bibliographystyle{ams-plain}

\bibliography{EndomorphismAlgebrasKugaSatake}

\end{document}